\documentclass[10pt,leqno]{amsart}
\usepackage{amsmath,amsthm,amsfonts,amsbsy,multicol}

\input{epsf.tex}
\usepackage{graphicx}
\usepackage{epstopdf}
\usepackage{color}
\DeclareGraphicsExtensions{.pdf, .jpg, .pnp, .bmp}
\usepackage{amssymb,latexsym}

\newtheorem{lemma}{Lemma}[section]

\newtheorem{theorem}[lemma]{Theorem}
\newtheorem{proposition}[lemma]{Proposition}
\newtheorem{definition}[lemma]{Definition}
\newtheorem{corollary}[lemma]{Corollary}
\newtheorem{example}[lemma]{Example}
\newtheorem{exercise}[lemma]{Exercise}
\newtheorem{remark}[lemma]{Remark}
\newtheorem{fig}[lemma]{Figure}
\newtheorem{tab}[lemma]{Table}

\newcommand{\bth}{\begin{theorem}}   \newcommand{\ethe}{\end{theorem}}
\newcommand{\bre}{\begin{remark}\em }   \newcommand{\ere}{\end{remark}}
\newcommand{\ble}{\begin{lemma}}      \newcommand{\ele}{\end{lemma}}
\newcommand{\bde}{\begin{definition}}   \newcommand{\ede}{\end{definition}}
\newcommand{\bco}{\begin{corollary}}     \newcommand{\eco}{\end{corollary}}
\newcommand{\bpr}{\begin{proposition}}  \newcommand{\epr}{\end{proposition}}
\newcommand{\bexer}{\begin{exercise}}     \newcommand{\eexer}{\end{exercise}}
\newcommand{\bexam}{\begin{example}}    \newcommand{\eexam}{\end{example}}
\newcommand{\bfi}{\begin{fig}}              \newcommand{\efi}{\end{fig}}
\newcommand{\btab}{\begin{tab}}         \newcommand{\etab}{\end{tab}}
\newcommand{\bpf}{\begin{proof}}        \newcommand{\epf}{\end{proof}}

\newcommand{\barr}{\begin{array}}   \newcommand{\earr}{\end{array}}
\newcommand{\beao}{\begin{eqnarray*}}       \newcommand{\eeao}{\end{eqnarray*}\noindent}
\newcommand{\beam}{\begin{eqnarray}}    \newcommand{\eeam}{\end{eqnarray}\noindent}
\newcommand{\beqq}{\begin{equation}}    \newcommand{\eeqq}{\end{equation}\noindent}

\newcommand{\bce}{\begin{center}}   \newcommand{\ece}{\end{center}}

 \newcommand{\un}{\underbrace}
\newcommand{\wt}{\widetilde}
\newcommand{\wh}{\widehat}
     
\newcommand{\mto}{m\to\infty}\newcommand{\nto}{n\to\infty}

\newcommand{\Rto}{R\to\infty}

\newcommand{\del}{\delta}
\newcommand{\D}{\Delta}
  \newcommand{\ep}{\epsilon}
\newcommand{\ka}{\kappa}
\newcommand{\lam}{\lambda}

\newcommand{\bbc}{{\mathcal C}} 

\newcommand{\bfE}{{\mathbb E}}\newcommand{\bbE}{{\mathcal E}}  
\newcommand{\bbf}{{\mathcal F}}

\newcommand{\bbi}{{\mathbb I}}

\newcommand{\bbl}{{\mathcal L}}
\newcommand{\bbm}{{\mathcal M}}
 \newcommand{\bbN}{{\mathbb N}}
\newcommand{\bfP}{{\mathbb P}}
                              
 \newcommand{\bbR}{{\mathbb R}}

\begin{document}

\bibliographystyle{plain}
\title[Approximating explicitly the mean reverting CEV process]{Approximating explicitly the mean reverting CEV process}

\author[N. Halidias]{N. Halidias}
\author[I. S. Stamatiou]{I. S. Stamatiou}
\address{Department of Mathematics,
University of the Aegean\\ Karlovassi, GR-83 200 Samos, Greece\\ Tel.  +3022730-82321, +3022730-82343}
\email{nick@aegean.gr, istamatiou@aegean.gr, joniou@gmail.com}

\begin{abstract}
In this paper we want to exploit further the semi-discrete method appeared in Halidias and Stamatiou (2015). We are interested in the numerical solution of mean reverting CEV processes that appear in financial mathematics models and are described as non negative solutions of certain stochastic differential equations with sub-linear diffusion coefficients of the form $(x_t)^q,$ where $\frac{1}{2}<q<1.$  Our goal is to construct explicit numerical schemes that preserve positivity. We  prove convergence of the proposed SD scheme with rate depending on the parameter $q.$ Furthermore, we verify our findings through numerical experiments and compare with other positivity preserving schemes. Finally, we show how to treat the whole two-dimensional stochastic volatility model, with instantaneous variance process given by the above mean reverting CEV process.
\end{abstract}
\date\today

\keywords{Explicit numerical scheme, mean reverting CEV process, positivity preserving, strong approximation error, order of convergence, stochastic volatility model. \newline{\bf AMS subject classification 2010:} 60H10, 60H35.}
\maketitle

\section{Introduction.}\label{sec:s1}
\setcounter{equation}{0}

Consider the following stochastic models
\beqq\label{eq1}
\left\{ \barr{ll}
S_t=S_0 + \int_{0}^{t}\mu\cdot S_udu + \int_{0}^{t}(V_u)^p \cdot S_u dW_u, & t\in [0,T],\\
&\\
V_t=V_0 + \int_0^t (k_1-k_2 V_s)ds + \int_{0}^{t}k_3(V_s)^q d\wt{W_s} & t\in [0,T],
\earr \right.
\eeqq
where $S_t$ represents the underlying financially observable variable, $V_t$ is the instantaneous volatility when $p=1$ or the instantaneous variance when $p=1/2$ and the Wiener  processes $W_t, \wt{W_t}$ have correlation $\rho.$

We assume that $V_t$ is a mean reverting CEV process of the above form, with the coefficients $k_i>0$ for $i=1,2,3$ and $q>1/2,$ since the process $V_t$ has to be non-negative.
To be more precise the above restriction on $q$ implies that $V_t$ is positive, i.e. $0$ is unattainable, as well as non explosive, i.e. $\infty$ is unattainable, as can be verified by the Feller's classification of boundaries \cite[Prop. 5.22]{karatzas_shreve:1988}. The steady-state level of $V_t$ is $k_1/k_2$ and the rate of mean reversion is $k_2$.

The system (\ref{eq1}) for $p=q=1/2$ is the Heston model. When $q=1$ we get the Brennan-Schwartz model \cite[Sec. II]{brennan_schwartz:1980}, which apparent its simple form, cannot provide analytical expressions for $S_t.$

Process $V_t$ for $q=1/2,$ also know as the CIR process \cite[Rel 13]{cox_et_al:1985}, by the initials of the authors that proposed it for the term structure of interest rates, has received a lot of attention and we just mention two latest contributions to the study of such processes (see \cite{alfonsi:2013}, \cite{halidias:2015b} and references therein).

Process $V_t$ for $1/2\leq q\leq1$ has been also considered for the dynamics of the short-term interest rate \cite[Rel (1)]{chan_et_al:1992}. The stationary distribution of the process has also been derived in \cite[Prop 2.2]{andersen_piterbarg:2007}.

We aim for a positive preserving scheme for the process $V_t.$ The scheme which we propose, and denote it semi-discrete (SD), preserves this analytical property of $V_t$ staying positive. The explicit Euler scheme fails to preserve positivity, as well as the standard Milstein scheme. We intend to apply the semi-discrete method for the numerical approximation of $V_t$ in model (\ref{eq1}) with $1/2<q<1$ and compare with other positivity preserving methods such as the balanced implicit method (BIM) (introduced by \cite[Rel (3.2)]{milstein_et_al:1998} with the positive preserving property \cite[Sec. 5]{kahl_schurz:2006}) and the balanced Milstein method (BMM) \cite[Th. 5.9]{kahl_schurz:2006}.\footnote{We give in the Appendix the form of all the above schemes for the approximation of $V_t$.} Finally, we  approximate the stochastic volatility model of (\ref{eq1}) with $p=1/2.$ In \cite{kahl_jackel:2006} a thorough treatment can be found, where also another stochastic volatility model is suggested.

Section \ref{sec:s2} provides the setting and the main results, Theorems \ref{t1} and \ref{t2}, concerning the $\bbl^2-$convergence of the proposed Semi Discrete (SD) method to the true solution of mean reverting CEV processes of the form of the stochastic volatility in (\ref{eq1}), as well as Theorem \ref{t3} which concerns the analogues of Theorems \ref{t1} and \ref{t2}, following an alternative approach. The main ingredient of this approach, inspired by \cite{halidias:2015b}, is the simplification of the numerical scheme proposed, by  altering the initial Brownian motion $(W_t)$ to another Brownian motion $(\hat{W}_t)$ justified by Levy's martingale characterization of Brownian motion, yielding to the same logarithmic rate as in Theorem \ref{t1}, but to a better polynomial rate $\frac{1}{2}(q-\frac{1}{2})$ instead of $\frac{1}{2}(q-\frac{1}{2})\wedge\frac{1}{8}$ as shown in  Theorem \ref{t2}.  

Section \ref{sec:s3} is devoted to the proof of Theorem \ref{t1}, while Section \ref{sec:s4} and \ref{sec:s5} concern the proofs of Theorems \ref{t2} and \ref{t3} respectively. Finally, Section \ref{sec:s6} presents illustrative figures where the behavior of the proposed scheme, regarding the order of convergence, is shown and a comparison with BIM and BMM schemes is given. In Section \ref{sec:s7} we treat the full model (\ref{eq1}) for a special case. Concluding remarks are in  Section \ref{sec:s8} and in Appendix \ref{A1} we briefly present  numerical schemes for the integration of the variance-volatility process $(V_t).$

\section{The setting and the main results.}\label{sec:s2}
\setcounter{equation}{0}

We consider the following SDE
\beqq\label{eq2}
x_t=x_0 + \int_0^t (k_1-k_2 x_s)ds + \int_{0}^{t}k_3(x_s)^q dW_s, \quad t\in [0,T],
\eeqq
where $k_1, k_2, k_3$ are positive and $1/2<q<1.$ Then, Feller's test implies that there is a unique strong solution such that $x_t>0$ a.s. when $x_0>0$ a.s. Let
\beqq\label{eq28}
f_{\theta}(x,y)= \underbrace{k_1 - k_2(1-\theta)x - \frac{k_3^2}{4(1+k_2\theta\D)}x^{2q-1}-k_2\theta y}_{f_1(x,y)} + \underbrace{\frac{k_3^2}{4(1+k_2\theta\D)}x^{2q-1}}_{f_2(x)}
\eeqq
and
\beqq\label{eq29}
g(x,y)=k_3 x^{q-\frac{1}{2}}\sqrt{y},
\eeqq
where $f(x,x)=a(x)=k_1-k_2 x$ and $g(x,x)=b(x)=k_3x^q.$

Let the partition $0=t_0<t_1<\ldots<t_N=T$ with $\D=T/N$ and consider the following process
$$
y_t^{SD}(q)=y_{t_n} + f_1(y_{t_n},y_t)\cdot\D + \int_{t_n}^tf_2(y_{t_n})ds + \int_{t_n}^{t} sgn(z_s)g(y_{t_n},y_s)dW_s,
$$
with $y_0=x_0$ a.s. or more explicitly
\beam
\nonumber
&&y_t^{SD}(q)=y_{t_n} + \left(k_1 - k_2(1-\theta)y_{t_n} - \frac{k_3^2}{4(1+k_2\theta\D)}(y_{t_n})^{2q-1}-k_2\theta y_t\right)\cdot\D\\
\label{eq3}&& + \int_{t_n}^{t}\frac{k_3^2}{4(1+k_2\theta\D)}(y_{t_n})^{2q-1} ds
+ k_3(y_{t_n})^{q-\frac{1}{2}}\int_{t_n}^{t} sgn(z_s)\sqrt{y_s}dW_s,
\eeam
for $ t\in(t_n,t_{n+1}],$ where $\theta\in[0,1]$ represents the level of implicitness and
\beqq\label{eq5}
z_t= \sqrt{y_n} +\frac{k_3}{2(1+k_2\theta\D)}(y_{t_n})^{q -\frac{1}{2}}(W_t-W_{t_n}),
\eeqq
with
\beqq\label{eq19}
y_n:=y_{t_n}\left(1- \frac{k_2\D}{1+k_2\theta\D}\right) + \frac{k_1\D}{1+k_2\theta\D} - \frac{k_3^2}{4(1+k_2\theta\D)^2}(y_{t_n})^{2q-1}\D.
\eeqq

Process (\ref{eq3}) is well defined when $y_n\geq0$ and this is true when $\frac{1}{(1+k_2\theta\D)}k_3^2\leq4(k_2\wedge k_1)$
and $\D(2-\theta)\leq\frac{1}{k_2}.$ Furthermore, (\ref{eq3}) has jumps at nodes $t_n.$ Solving for $y_t,$ we end up with the following explicit scheme
\beqq\label{eq4}
y_t^{SD}(q)=y_n + \int_{t_n}^{t}\frac{k_3^2}{4(1+k_2\theta\D)^2}(y_{t_n})^{2q-1}ds
+ \frac{k_3}{1+k_2\theta\D}(y_{t_n})^{q-\frac{1}{2}}\int_{t_n}^{t} sgn(z_s)\sqrt{y_s} dW_s,
\eeqq
with solution in each step given by \cite[Rel. (4.39), p.123]{kloeden_platen:1995}
$$
y_t^{SD}(q)=(z_t)^{2},
$$
which has the pleasant feature $y_t^{SD}(q)\geq0.$

\fbox{\begin{minipage}[t] {15.8cm}
\textbf{Assumption A}  Let the parameters $k_1,k_2,k_3$ be positive and such that $\frac{1}{(1+k_2\theta\D)}k_3^2\leq4(k_2\wedge k_1)$ and
 consider $\D>0$ such that $\D(2-\theta)<\frac{1}{k_2},$ for $\theta\in[0,1].$ Moreover assume $x_0\geq0$ a.s. and $\bfE(x_0)^p<A$ for some $p\geq4.$
\end{minipage}}

\bth \label{t1}[Logarithmic rate of convergence]
Let Assumption A hold. The semi-discrete scheme (\ref{eq4}) converges to the true solution of (\ref{eq2}) in the mean square sense with rate given by
\beqq \label{eq6}
\bfE\sup_{0\leq t\leq T}|y_t^{SD}(q)-x_t|^2\leq \frac{C}{\sqrt{\ln (\D)^{-1}}},
\eeqq
where $C$ is independent of $\D$ and particularly
$$C:=72\sqrt{\frac{6}{\ep}}(k_3)^4T^2e^{6T^2k_2+k_2T},$$
where $\ep$ is such that $0<\ep<\frac{1}{4}\wedge (q-\frac{1}{2}).$
\ethe

\fbox{\begin{minipage}[t] {15.8cm}
\textbf{Assumption B}  Let Assumption A hold where now  $x_0\in\bbR$ and $x_0>0.$
\end{minipage}}

\bth \label{t2}[Polynomial rate of convergence]
Let Assumption B hold. The semi-discrete scheme (\ref{eq4}) converges to the true solution of (\ref{eq2}) in the mean square sense with rate given by,
\beqq \label{eq30}
\bfE\sup_{0\leq t\leq T}|y_t^{SD}(q)-x_t|^2\leq C\D^{(q-\frac{1}{2})\wedge\frac{1}{4}},
\eeqq
where
\beao
C&:=&8\left(12k_3^2T(\sqrt{A_{4q}(x_0+ k_1T)^{4q}C_{k_2,k_3,\theta,\D}}\vee1) \bigvee 3k_3^2T\sqrt{A_2(x_0 + k_1T)^2}\sqrt{\hat{A}_{4q-2}}\right)\\
&&\times\left(2e^{6k_2T^2} + \frac{C_{HK}}{\ep-1}(x_0)^{(1-q)\nu(\lam)}\right),
\eeao
and $C_{HK}$ is a constant described  in (\ref{eq47}) and $\lam$ is an appropriately chosen positive parameter which satisfies (\ref{eq43}) and always exist, $\nu(\lam):=\frac{\lam}{2(1-q)^2(k_3)^2} - 1,$ quantity $C_{k_2,k_3,\theta,\D}$ is given in Lemma \ref{l1} and $\ep>1.$
\ethe

Inspired by  \cite{halidias:2015b} we remove the term $sgn(z_s)$ from (\ref{eq3}) by considering the process
\beqq\label{eq101}
 \wt{W}_t=\int_{0}^{t} sgn(z_s)dW_s,
\eeqq
which is a martingale with quadratic variation $< \wt{W}_t, \wt{W}_t>=t$ and thus a standard Brownian motion w.r.t. its own filtration, justified by Levy's theorem \cite[Th. 3.16, p.157]{karatzas_shreve:1988}. Therefore, the compact form of (\ref{eq3}) becomes
\beao
&&y_t^{SD}=x_0 + \int_0^t\left(k_1 - k_2(1-\theta)y_{\hat{s}}-k_2\theta y_{\wt{s}}\right)ds\\
&& + \int_t^{t_{n+1}}\left(k_1 - k_2(1-\theta)y_{t_n} - \frac{k_3^2}{4(1+k_2\theta\D)}(y_{t_n})^{2q-1}-k_2\theta y_{t}\right)ds + k_3\int_{0}^{t} s(y_{\hat{s}})^{q-\frac{1}{2}}\sqrt{y_s}d\wt{W}_s,
\eeao
for $t\in(t_n,t_{n+1}].$ Consider also the process
\beqq\label{eq102}
\wt{x}_t=x_0 + \int_0^t (k_1-k_2 \wt{x}_s)ds + \int_{0}^{t}k_3(\wt{x}_s)^q d\wt{W}_s, \quad t\in [0,T].
\eeqq
The process $(x_t)$ of (\ref{eq2}) and the process $(\wt{x}_t)$ of (\ref{eq102}) have the same distribution. We show in the following that $\bfE\sup_{0\leq t\leq T}|y_t^{SD}(q)-\wt{x}_t|^2\rightarrow0$ as $\D\downarrow0$ thus the same holds for the unique solution of (\ref{eq2}), i.e. $\bfE\sup_{0\leq t\leq T}|y_t^{SD}(q)-x_t|^2\rightarrow0$ as $\D\downarrow0.$ To simplify notation we write $\wt{W}, (\wt{x}_t)$ as $W, (x_t).$
We end up with the following explicit scheme
\beqq\label{eq103}
y_t^{SD}(q)=y_n + \int_{t_n}^{t}\frac{k_3^2}{4(1+k_2\theta\D)^2}(y_{t_n})^{2q-1}ds
+ \frac{k_3}{1+k_2\theta\D}(y_{t_n})^{q-\frac{1}{2}}\int_{t_n}^{t} \sqrt{y_s} dW_s,
\eeqq
where $y_n$ is as in (\ref{eq19}).

\bth \label{t3}[Logarithmic and Polynomial rate of convergence]

Let Assumption A hold. The semi-discrete scheme (\ref{eq103}) converges to the true solution of (\ref{eq2}) in the mean square sense with rate given by
\beqq \label{eq104}
\bfE\sup_{0\leq t\leq T}|y_t^{SD}(q)-x_t|^2\leq \frac{C}{\sqrt{\ln (\D)^{-1}}},
\eeqq
where $C$ is independent of $\D$ and given by
$$C:=32\sqrt{\frac{6}{\ep}}(k_3)^4T^2e^{6T^2k_2+k_2T},$$
where $\ep$ is such that $0<\ep<q-\frac{1}{2}.$

In case  Assumption B holds, the semi-discrete scheme (\ref{eq103}) converges to the true solution of (\ref{eq2}) in the mean square sense with rate given by,
\beqq \label{eq105}
\bfE\sup_{0\leq t\leq T}|y_t^{SD}(q)-x_t|^2\leq C\D^{(q-\frac{1}{2})},
\eeqq
where
\beao
C&:=&\frac{16}{3}\left(12k_3^2T(\sqrt{A_{4q}(x_0+ k_1T)^{4q}C_{k_2,k_3,\theta,\D}}\vee1) \bigvee 3k_3^2T\sqrt{A_2(x_0 + k_1T)^2}\sqrt{\hat{A}_{4q-2}}\right)\\
&&\times\left(2e^{6k_2T^2} + \frac{C_{HK}}{\ep-1}(x_0)^{(1-q)\nu(\lam)}\right),
\eeao
and $C_{HK}$ is the constant described  in (\ref{eq47}) and $\lam$ is an appropriately chosen positive parameter which satisfies (\ref{eq43}) and always exist, $\nu(\lam):=\frac{\lam}{2(1-q)^2(k_3)^2} - 1,$ quantity $C_{k_2,k_3,\theta,\D}$ is given in Lemma \ref{l1} and $\ep>1.$
\ethe
In the following sections we write for simplicity $y_t^{SD}$ or $y_t$ for $y_t^{SD}(q).$

\section{Logarithmic rate of convergence.}\label{sec:s3}
\setcounter{equation}{0}

We rewrite (\ref{eq3}) in a compact form
\beao
&&y_t^{SD}=x_0 + \int_0^t\left(k_1 - k_2(1-\theta)y_{\hat{s}}-k_2\theta y_{\wt{s}}\right)ds\\
&& + \int_t^{t_{n+1}}\left(k_1 - k_2(1-\theta)y_{t_n} - \frac{k_3^2}{4(1+k_2\theta\D)}(y_{t_n})^{2q-1}-k_2\theta y_{t}\right)ds + k_3\int_{0}^{t} sgn(z_s)(y_{\hat{s}})^{q-\frac{1}{2}}\sqrt{y_s}dW_s,
\eeao
for $t\in(t_n,t_{n+1}]$ where
$$
\hat{s}=t_j, s\in(t_j,t_{j+1}], \, j=0,\ldots,n, \quad
\wt{s}=\left\{ \barr{ll}  t_{j+1},  & \mbox{for } \, s\in[t_j,t_{j+1}], \\  t,  & \mbox{for }\,  s\in[t_n,t] \earr j=0,\ldots,n-1.\right.
$$

\subsection{Moment bounds}\label{ssec:s3.1}
\setcounter{equation}{0}

\ble[Moment bound for SD approximation]\label{l4}
It holds that
$$
\bfE\sup_{0\leq t\leq T}(y_t)^p\leq A_p\bfE(x_0 + k_1T)^p,
$$
for any $p>2,$ where $A_p:= \exp\left\{\frac{p(p-1)}{2}k_3^2\left(\frac{p-1}{2p} + \frac{2^{p-1}}{p}\right)T\right\}.$
\ele
\bpf[Proof of Lemma \ref{l4}]
We first observe that $(y_t)$ is bounded in the following way
\beao
0\leq y_t&\leq& x_0 + \int_0^t k_1ds + k_3\int_{0}^{t} sgn(z_s)(y_{\hat{s}})^{q-\frac{1}{2}}\sqrt{y_s}dW_s\\
&\leq& x_0 + k_1T + k_3\int_{0}^{t} sgn(z_s)(y_{\hat{s}})^{q-\frac{1}{2}}\sqrt{y_s}dW_s:=u_t
\eeao
a.s., where the lower bound comes from the construction of $(y_t)$ and the upper bound follows from a comparison theorem. We will bound $(u_t)$ and therefore $(y_t),$ since $0\leq y_t\leq u_t$ a.s. Set the stopping time $\tau_R:=\inf\{t\in[0,T]: u_t>R\},$ for $R>0$ with the convention $\inf \emptyset=\infty.$ Application of Ito's formula on $(u_{t\wedge\tau_R})^p$ implies
\beao
(u_{t\wedge\tau_R})^p&=&(x_0 + k_1T)^p + \frac{p(p-1)}{2}k_3^2\int_{0}^{t\wedge\tau_R} (u_s)^{p-2}(y_{\hat{s}})^{2q-1}y_sds\\
&& + pk_3\int_{0}^{t\wedge\tau_R}sgn(z_s)(u_s)^{p-1}(y_{\hat{s}})^{q-\frac{1}{2}}\sqrt{y_s}dW_s\\
&\leq&(x_0 + k_1T)^p + \frac{p(p-1)}{2}k_3^2\int_{0}^{t\wedge\tau_R} (u_s)^{p-1}(y_{\hat{s}})^{2q-1}ds +
M_t\\
&\leq&(x_0 + k_1T)^p + \frac{p(p-1)}{2}k_3^2\int_{0}^{t\wedge\tau_R} \left(\frac{p-1}{2p}(u_s)^{p} + \frac{2^{p-1}}{p}(y_{\hat{s}})^{(2q-1)p}\right)ds +  M_t\\
&\leq&(x_0 + k_1T)^p + \frac{p(p-1)}{2}k_3^2\left(\frac{p-1}{2p} + \frac{2^{p-1}}{p}\right)\int_{0}^{t\wedge\tau_R} (u_s)^{p}ds +  M_t,
\eeao
where in the second step we have used that $0\leq y_t\leq u_t,$ in the third step the inequality $x^{p-1}y\leq \ep\frac{p-1}{p}x^p + \frac{1}{p\ep^{p-1}}y^p,$ valid for $x\wedge y\geq0$ and $p>1$ with $\ep=\frac{1}{2},$ in the final step the fact $\frac{1}{2}<q<1$ and $M_t:=pk_3\int_{0}^{t\wedge\tau_R}sgn(z_s)(u_s)^{p-1}(y_{\hat{s}})^{q-\frac{1}{2}}\sqrt{y_s}dW_s.$ Taking expectations in the above inequality and using that $M_t$ is a local martingale vanishing at $0,$ we get
\beao
\bfE(u_{t\wedge\tau_R})^p&\leq&\bfE(x_0 + k_1T)^p + \frac{p(p-1)}{2}k_3^2\left(\frac{p-1}{2p} + \frac{2^{p-1}}{p}\right)\int_{0}^{t} \bfE(u_{s\wedge\tau_R})^{p}ds\\
&\leq&\bfE(x_0 + k_1T)^p \exp\left\{\frac{p(p-1)}{2}k_3^2\left(\frac{p-1}{2p} + \frac{2^{p-1}}{p}\right)T\right\}\\
&\leq&A_p\bfE(x_0 + k_1T)^p,
\eeao
where we have applied Gronwall inequality \cite[Rel. 7]{gronwall:1919}. We have that
$$
(y_{t\wedge\tau_R})^p=(y_{\tau_R})^p\bbi_{\{\tau_R\leq t\}} + (y_{t})^p\bbi_{\{t<\tau_R\}}\geq(y_{t})^p\bbi_{\{t<\tau_R\}},
$$
thus taking expectations in the above inequality and using the estimated upper bound for $\bfE(u_{t\wedge\tau_R})^p$ we arrive at
$$
\bfE(y_{t})^p\bbi_{\{t<\tau_R\}}\leq \bfE(y_{t\wedge\tau_R})^p\leq \bfE(u_{t\wedge\tau_R})^p\leq A_p\bfE(x_0 + k_1T)^p,
$$
and taking the limit as $\Rto,$ we get
$$
\lim_{\Rto}\bfE(y_{t})^p\bbi_{\{t<\tau_R\}}\leq A_p\bfE(x_0 + k_1T)^p.
$$
Let us fix $t.$ The sequence of stopping times $\tau_R$ is increasing in $R$ and $t\wedge\tau_R\rightarrow t$ as $\Rto,$ thus the sequence $(y_{t})^p\bbi_{\{t<\tau_R\}}$ is nondecreasing in $R$ and $(y_{t})^p\bbi_{\{t<\tau_R\}}\rightarrow(y_{t})^p$ as $\Rto.$ Application of the monotone convergence theorem implies
\beqq\label{eq13}
\bfE(y_{t})^p\leq A_p\bfE(x_0 + k_1T)^p,
\eeqq
for any $p>2.$ Using again Ito's formula on $(u_t)^p$, taking the supremum and then using Doob's martingale inequality on the diffusion term we bound $\bfE\sup_{0\leq t\leq T}(u_t)^p$ and thus $\bfE\sup_{0\leq t\leq T}(y_t)^p.$
\epf

\ble[Error bound for SD scheme]\label{l5}
Let $n_s$ integer such that $s\in[t_{n_s},t_{n_s+1}].$ It holds that
$$
\bfE|y_s-y_{\hat{s}}|^p\leq\hat{A}_p\D^{p/2}, \qquad \bfE |y_s-y_{\wt{s}}|^p<\wt{A}_p\D^{p/2},
$$
for any $p>0,$ where the positive quantities $\hat{A}_p,\wt{A}_p$ do not depend on $\D.$
\ele
\bpf[Proof of Lemma \ref{l5}]
First we take a $p>2.$ It holds that
\beao
&&|y_s-y_{\hat{s}}|^p=\Big|\int_{t_{\wh{n_s}}}^s\left(k_1 - k_2(1-\theta)y_{\hat{u}}-k_2\theta y_{\wt{u}}\right)du+ \int_{t_{\wh{n_s}}}^{t_{n_s+1}}k_2\theta y_{\hat{s}}du  -\int_s^{t_{n_s+1}}k_2\theta y_{s}du\\
&&+ \int_s^{t_{\wh{n_s}}}\left(k_1 - k_2(1-\theta)y_{t_{n_s}} - \frac{k_3^2}{4(1+k_2\theta\D)}(y_{t_{n_s}})^{2q-1}\right)du + k_3\int_{t_{\wh{n_s}}}^{s} sgn(z_u)(y_{\hat{u}})^{q-\frac{1}{2}}\sqrt{y_u}dW_u\Big|^p\\
&\leq&5^{p-1}\Big(\big|\int_{t_{\wh{n_s}}}^s\left(k_1 - k_2(1-\theta)y_{\hat{u}}-k_2\theta y_{\wt{u}}\right)du\big|^p + k_2^p\theta^p(y_{\hat{s}})^p(t_{n_s+1}-t_{\wh{n_s}})^p + k_2^p\theta^p(y_{s})^p(t_{n_s+1}-s)^p\\
&&+ \left|\int_s^{t_{\wh{n_s}}}\left(k_1 - k_2(1-\theta)y_{t_{n_s}} - \frac{k_3^2}{4(1+k_2\theta\D)}(y_{t_{n_s}})^{2q-1}\right)du\right|^p + k_3^p\big|\int_{t_{\wh{n_s}}}^{s} sgn(z_u)(y_{\hat{u}})^{q-\frac{1}{2}}\sqrt{y_u}dW_u\big|^p\Big)\\
&\leq&5^{p-1}\Big(|t_{\wh{n_s}}-s|^{p-1}\int_{t_{\wh{n_s}}}^s\left|k_1 - k_2(1-\theta)y_{\hat{u}}-k_2\theta y_{\wt{u}}\right|^pdu + k_2^p\theta^p\left((y_{\hat{s}})^p+(y_{s})^p\right)\D^p\\
&&+ \left|k_1 - k_2(1-\theta)y_{t_{n_s}} - \frac{k_3^2}{4(1+k_2\theta\D)}(y_{t_{n_s}})^{2q-1}\right|^p\D^p + k_3^p\big|\int_{t_{\wh{n_s}}}^{s} sgn(z_u)(y_{\hat{u}})^{q-\frac{1}{2}}\sqrt{y_u}dW_u\big|^p\Big),
\eeao
where we have used Cauchy-Schwarz inequality. Taking expectations in the above inequality and using Lemma \ref{l4} and Doob's martingale inequality on the diffusion term we conclude
\beqq\label{eq14}
\bfE|y_s-y_{\hat{s}}|^p\leq\hat{A}_p\D^{p/2},
\eeqq
where the positive quantity $\hat{A}_p$ except on $p,$ depends also on the parameters $k_1,k_2,k_3,\theta,q,$ but not on $\D.$ Now, for  $0<p<2$ we get
$$
\bfE|y_s-y_{\hat{s}}|^p\leq\left(\bfE|y_s-y_{\hat{s}}|^2\right)^{p/2}\leq\hat{A}_p\D^{p/2},
$$
where we have used Jensen inequality for the concave function $\phi(x)=x^{p/2}.$ Following the same lines, we can show that
\beqq\label{eq15}
\bfE|y_s-y_{\wt{s}}|^p\leq\wt{A}_p\D^{p/2},
\eeqq
for any  $0<p,$ where the positive quantity $\wt{A}_p$ except on $p,$ depends also on the parameters $k_1,k_2,k_3,\theta,q,$ but not on $\D.$
\epf

For the rest of this section we rewrite again the compact form of (\ref{eq3}) in the following way
\beqq\label{eq16}
y_t^{SD}=\un{x_0 + \int_0^t f_{\theta}(y_{\hat{s}},y_{\wt{s}})ds + \int_{0}^{t}sgn(z_s)g(y_{\hat{s}},y_s)dW_s}_{h_t} + \int_t^{t_{n+1}}f_1(y_{t_n},y_{t})ds,
\eeqq
where $f_{\theta}(\cdot,\cdot)$ is given by (\ref{eq28}) and the auxiliary process $(h_t)$ is close to $(y_t)$ as shown in the next result.

\ble[Moment bounds involving the auxiliary process]\label{l6}
For any $s\in[0,T]$ it holds that
\beqq\label{eq17}
\bfE|h_s-y_{s}|^p\leq C_p\D^{p},  \quad \bfE|h_s|^p\leq C_h,
\eeqq
and for $s\in[t_{n},t_{n+1}]$ we have that
\beqq\label{eq18}
\bfE|h_s-y_{\hat{s}}|^p\leq \hat{C}_p\D^{p/2}, \qquad \bfE|h_s-y_{\wt{s}}|^p\leq \wt{C}_p\D^{p/2},
\eeqq
for any $p>0,$ where  the positive quantities $C_p,\hat{C}_p,\wt{C}_p,C_h$ do not depend on $\D.$
\ele
\bpf[Proof of Lemma \ref{l6}]
We have that
$$
|h_s-y_s|^p=\left|\int_s^{t_{n+1}}f_1(y_{t_n},y_{t})du \right|^p\leq|t_{n+1}-s|^p|f_1(y_{t_n},y_{t})|^p,
$$
for any $p>0,$ where we have used (\ref{eq16}). Using Lemma \ref{l4} we get the left part of (\ref{eq17}). Now for $p>2$ and noting that
\beao
\bfE|h_s|^p&\leq&2^{p-1}\bfE|h_s-y_s|^p + 2^{p-1}\bfE|y_s|^p\\
&\leq&2^{p-1}C_p\D^{p} + 2^{p-1}A_p\bfE(x_0 + k_1T)^p\leq C_h,
\eeao
we get the right part of (\ref{eq17}), where we have used Lemma \ref{l4}. The case $0<p<2$ follows by Jensen's inequality as in Lemma \ref{l5}.

Furthermore, for $s\in[t_{n},t_{n+1}]$ and $p>2$ it holds
\beao
\bfE|h_s-y_{\hat{s}}|^p&\leq&2^{p-1}\bfE|h_s-y_s|^p + 2^{p-1}\bfE|y_s-y_{\hat{s}}|^p\\
&\leq& 2^{p-1}C_p\D^{p} + 2^{p-1}\hat{A}_p\D^{p/2}\leq\hat{C}_p\D^{p/2}
\eeao
where we have used (\ref{eq14}) and in the same manner
$$
\bfE|h_s-y_{\wt{s}}|^p\leq 2^{p-1}C_p\D^{p} + 2^{p-1}\wt{A}_p\D^{p/2}\leq\wt{C}_p\D^{p/2}.$$
The case $0<p<2$ follows by Jensen's inequality.
\epf

\subsection{Convergence of the auxiliary process $(h_t)$ to $(x_t)$ in $\bbl^1$}\label{ssec:s3.2}

We first estimate the probability of $z_t$ being negative when at the same time $y_{t_n}>\D^{1-2\xi},$ for $0<\xi<\frac{1}{2}.$
\ble\label{l1}
For every $t\in[t_n,t_{n+1}]$ it holds
\beqq\label{eq7}
\bfP\left(\{z_t\leq0\}\cap\{y_{t_n}>\D^{1-2\xi}\}\right)\leq C_{k_2,k_3,\theta,\D}\sqrt{\D},
\eeqq
where $C_{k_2,k_3,\theta,\D}:=\frac{k_3}{\sqrt{1- k_2(2-\theta)\D}}$ and $\D(2-\theta)<\frac{1}{k_2}$ and $\frac{k_3^2}{(1+k_2\theta\D)}\leq4k_2.$
Relation (\ref{eq7}) implies that $\bfP\left(\{z_t\leq0\}\cap\{y_{t_n}>\D^{1-2\xi}\}\right)=O(\sqrt{\D}),$ as $\D\downarrow0.$
\ele

\bpf[Proof of Lemma \ref{l1}]
By the definition (\ref{eq5}) of $(z_t)$ for $t\in[t_n,t_{n+1}]$ and for $0<\xi<\frac{1}{2},$ we have that
\beam
\nonumber A:=\{z_t\leq0\}\cap\{y_{t_n}>\D^{1-2\xi}\}&=&\left\{(y_{t_n})^{q -\frac{1}{2}}(W_t-W_{t_n})\leq -\frac{2(1+k_2\theta\D)}{k_3}\sqrt{y_n}\right\}\cap\{y_{t_n}>\D^{1-2\xi}\}\\
\label{eq8}&\subseteq&A_1 \cup A_2,
\eeam
where
$$
A_1:=\left\{W_t-W_{t_n}\leq -\frac{2(1+k_2\theta\D)}{k_3}\sqrt{y_n}(y_{t_n})^{-q + \frac{1}{2}}\right\}\cap\{y_{t_n}\geq1\},
$$
and
$$
A_2:=\left\{W_t-W_{t_n}\leq -\frac{2(1+k_2\theta\D)}{k_3}\sqrt{y_n}(y_{t_n})^{-q + \frac{1}{2}}\right\}\cap\{1>y_{t_n}>\D^{1-2\xi}\}.$$
The following inclusion relations hold for the event $A_1,$
\beao
&&A_1\\
&\subseteq&\left\{\D W_{n}\leq -\frac{2(1+k_2\theta\D)}{k_3}\sqrt{y_{t_n}\left(1- \frac{k_2\D}{1+k_2\theta\D}\right) - \frac{(k_3)^2\D}{4(1+k_2\theta\D)^2}(y_{t_n})^{2q-1}}(y_{t_n})^{-q + \frac{1}{2}}\right\}\cap\{y_{t_n}\geq1\}\\
&\subseteq&\left\{\D W_{n}\leq -\frac{2(1+k_2\theta\D)}{k_3}\sqrt{\frac{1- k_2(2-\theta)\D}{1+k_2\theta\D}- \frac{(k_3)^2\D}{4(1+k_2\theta\D)^2}}\right\}\\
&\subseteq&\left\{\frac{\D W_{n}}{\sqrt{t- t_n}}\leq -\frac{2}{k_3}\frac{\sqrt{(1- k_2(2-\theta)\D)(1+k_2\theta\D)}}{\sqrt{t- t_n}}\right\}
\eeao
when $\D(2-\theta)<\frac{1}{k_2}$ and $\frac{(k_3)^2}{(1+k_2\theta\D)}\leq4k_2,$ where $\D W_{n}:=W_t-W_{t_n}.$ It holds that
\beqq\label{eq8.1}
\bfP(G\leq -\beta)=\int_{-\infty}^{-\beta}\frac{1}{\sqrt{2\pi}}e^{-u^2/2}du\leq\int_{-\infty}^{-\beta}e^{-u^2/2}du=\int_{\beta}^{\infty}e^{-u^2/2}du\leq \frac{1}{\beta}e^{-(\beta)^2/2},
\eeqq
for every standard normal random variable $G$, where in the last step we have used \cite[Ineq. (9.20), p.112]{karatzas_shreve:1988} valid for $\beta>0$. Using the fact that $\frac{\D W_{n}}{\sqrt{t- t_n}}$ is  a standard normal r.v. and ignoring the exponential term in (\ref{eq8.1}), since its exponent is negative,  we get that
\beqq\label{eq8.2}
\bfP(A_1)\leq\frac{k_3}{2\sqrt{(1- k_2(2-\theta)\D)}}\sqrt{t- t_n}\leq C_{k_2,k_3,\theta,\D}\sqrt{\D}.
\eeqq
The following inclusion relations hold for the event $A_2,$
\beao
&&A_2\\
&\subseteq&\left\{\D W_{n}\leq -\frac{2(1+k_2\theta\D)}{k_3}\sqrt{y_{t_n}\frac{1- k_2(1-\theta)\D}{1+k_2\theta\D} + \frac{k_1\D}{1+k_2\theta\D} - \frac{(k_3)^2\D}{4(1+k_2\theta\D)^2}(y_{t_n})^{2q-1}}(y_{t_n})^{-q + \frac{1}{2}}\right\}\\
&&\cap\{1>y_{t_n}>\D^{1-2\xi}\}\\
&\subseteq&\left\{\D W_{n}\leq -\frac{2(1+k_2\theta\D)}{k_3}\sqrt{\D^{1-2\xi}\frac{1- k_2(1-\theta)\D}{1+k_2\theta\D} + \left(k_1- \frac{(k_3)^2}{4(1+k_2\theta\D)}\right)\frac{\D}{1+k_2\theta\D}}\right\}\\
&\subseteq&\left\{\frac{\D W_{n}}{\sqrt{t- t_n}}\leq -\frac{2}{k_3}\frac{\sqrt{(1- k_2(1-\theta)\D)(1+k_2\theta\D)}}{\sqrt{t- t_n}}\D^{\frac{1}{2}-\xi}\right\}
\eeao
when $\D(1-\theta)<\frac{1}{k_2}$ and $\frac{(k_3)^2}{(1+k_2\theta\D)}\leq4k_1.$  Using again (\ref{eq8.1}) we have that
\beam
\nonumber\bfP(A_2)&\leq& \frac{k_3}{2\sqrt{(1- k_2(1-\theta)\D)}}\sqrt{t-t_n}\D^{\xi-\frac{1}{2}}e^{-\frac{2}{(k_3)^2}\frac{(1- k_2(1-\theta)\D)(1+k_2\theta\D)}{\sqrt{t- t_n}}\D^{1-2\xi}}\\
\label{eq8.3}&\leq&\frac{k_3}{2\sqrt{(1- k_2(1-\theta)\D)}}\D^{\xi}e^{-\frac{2}{(k_3)^2}(1- k_2(1-\theta)\D)(1+k_2\theta\D)\D^{-2\xi}}.
\eeam
Taking probabilities in the inclusion relation (\ref{eq8}) and using (\ref{eq8.2}) and (\ref{eq8.3}) we get
\beao
\bfP(A)&\leq& \bfP(A_1) + \bfP(A_2)\\
&\leq& C_{k_2,k_3,\theta,\D}\sqrt{\D} + \frac{k_3}{2\sqrt{(1- k_2(1-\theta)\D)}}\D^{\xi}e^{-\frac{2}{(k_3)^2}(1- k_2(1-\theta)\D)(1+k_2\theta\D)\D^{-2\xi}}\\
&\leq& C_{k_2,k_3,\theta,\D}\sqrt{\D},
\eeao
since $\D^{\xi}e^{-\D^{-2\xi}}=o(\sqrt{\D})$ as $\D\downarrow0.$ Finally, note that $C_{k_2,k_3,\theta,\D}= \frac{k_3}{\sqrt{1- k_2(2-\theta)\D}}\rightarrow k_3$ as $\D\downarrow0$ which justifies the $O(\sqrt{\D})$ notation, (see for example \cite{olver:1997}).
\epf

We will use the representation (\ref{eq16}) and write
\beqq\label{eq20}
h_t-x_t=\int_0^t \left(f_{\theta}(y_{\hat{s}},y_{\wt{s}})-f_{\theta}(x_s,x_s)\right)ds + \int_{0}^{t}\left(sgn(z_s)g(y_{\hat{s}},y_s)-g(x_s,x_s)\right)dW_s.
\eeqq

\bpr\label{pr1}
Let Assumption A hold. Then we have
\beqq \label{eq21}
\sup_{0\leq t\leq T}\bfE|h_{t}-x_{t}|\leq\left(J_2\frac{\D^{\frac{1}{4}}}{me_m} + J_3\frac{\D^{q-\frac{1}{2}}}{me_m} +3 k_3^2T\frac{1}{m}\right)e^{k_2T},
\eeqq
for any $m>1,$ where $e_m=e^{-m(m+1)/2}$ and
$$J_2:=12k_3^2T(\sqrt{A_{4q}\bfE(x_0+ k_1T)^{4q}C_{k_2,k_3,\theta,\D}}\vee1), \,J_3:=3k_3^2T\sqrt{A_2\bfE(x_0 + k_1T)^2}\sqrt{\hat{A}_{4q-2}}.$$
\epr

\bpf[Proof of Proposition \ref{pr1}]
Let the non increasing sequence $\{e_m\}_{m\in\bbN}$ with $e_m=e^{-m(m+1)/2}$ and $e_0=1.$ We introduce the following sequence of smooth approximations of $|x|,$ (method of Yamada and Watanabe, \cite{yamada_watanabe:1971})
$$
\phi_m(x)=\int_0^{|x|}dy\int_0^{y}\psi_m(u)du,
$$
where the existence of the continuous function $\psi_m(u)$ with $0\leq \psi_m(u) \leq 2/(mu)$ and support in $(e_m,e_{m-1})$ is justified by $\int_{e_m}^{e_{m-1}}(du/u)=m.$ The following relations hold for $\phi_m\in\bbc^2(\bbR,\bbR)$ with $\phi_m(0)=0,$
 $$
 |x| - e_{m-1}\leq\phi_m(x)\leq |x|, \quad |\phi_{m}^{\prime}(x)|\leq1, \quad x\in\bbR, $$
 $$
 |\phi_{m}^{\prime \prime }(x)|\leq\frac{2}{m|x|}, \,\hbox{ when }  \,e_m<|x|<e_{m-1} \,\hbox{ and }  \,  |\phi_{m}^{\prime \prime }(x)|=0 \,\hbox{ otherwise.}
 $$
We have that
\beqq\label{eq22}
\bfE|h_{t}-x_{t}| \leq e_{m-1} + \bfE\phi_m(h_{t}-x_{t}).
\eeqq
It holds that
\beam\nonumber
f_{\theta}(y_{\hat{s}},y_{\wt{s}})-f_{\theta}(x_s,x_s)&=&(k_1 - k_2(1-\theta)y_{\hat{s}}-k_2\theta y_{\wt{s}}) - (k_1-k_2x_s)\\
\nonumber&=& - k_2(1-\theta)(y_{\hat{s}}-x_s)-k_2\theta (y_{\wt{s}}-x_s)\\
\label{eq23}&=&  k_2(1-\theta)(h_s - y_{\hat{s}}) + k_2\theta (h_s-y_{\wt{s}}) -k_2(h_s-x_s)
\eeam
and
\beam\nonumber
&&|sgn(z_s)g(y_{\hat{s}},y_s)-g(x_s,x_s)|^2=|sgn(z_s)k_3(y_{\hat{s}})^{q-\frac{1}{2}}\sqrt{y_s}-k_3(x_s)^q|^2\\
\nonumber&\leq&k_3^2\left( (y_{\hat{s}})^{q-\frac{1}{2}}\sqrt{y_s}(sgn(z_s)-1) + \sqrt{y_s}\left((y_{\hat{s}})^{q-\frac{1}{2}}-(y_{s})^{q-\frac{1}{2}}\right)   + ((y_s)^q-(x_s)^q)\right)^2\\
\nonumber&\leq&3 k_3^2 \left( (y_{\hat{s}})^{2q-1}y_s(sgn(z_s)-1)^2 + y_s\left((y_{\hat{s}})^{q-\frac{1}{2}}-(y_{s})^{q-\frac{1}{2}}\right)^2  + ((y_s)^q-(x_s)^q)^2\right)\\
\nonumber&\leq&3 k_3^2 \left( (y_{\hat{s}})^{2q-1}y_s(sgn(z_s)-1)^2 + y_s|y_{\hat{s}}-y_{s}|^{2q-1} + (\sqrt{|y_s-x_s|})^2\right)\\
\label{eq24}&\leq&3 k_3^2 \left( (y_{\hat{s}})^{2q-1}y_s(sgn(z_s)-1)^2 + y_s|y_{\hat{s}}-y_{s}|^{2q-1} + |h_s-y_s| + |h_s-x_s|\right),
\eeam
where we have used properties of Holder continuous functions and namely the fact that $x^q$ is $q-$Holder continuous for $q\leq1$, i.e. $|x^q-y^q|\leq|x-y|^q,$ and   that $x^q$ is $1/2-$Holder continuous since $q>1/2.$
Application of Ito's formula to the sequence $\{\phi_m\}_{m\in\bbN},$ implies
\beao
&&\phi_m(h_{t}-x_{t})= \int_{0}^{t} \phi_m^{\prime}(h_s-x_s) (f_{\theta}(y_{\hat{s}},y_{\wt{s}})-f_{\theta}(x_s,x_s))ds +M_t\\
&&+ \frac{1}{2}\int_{0}^{t} \phi_m^{\prime\prime}(h_s-x_s) (sgn(z_s)g(y_{\hat{s}},y_{\wt{s}})-g(x_s,x_s))^2 ds \\
&\leq& \int_{0}^{t} \left(k_2(1-\theta)|h_s - y_{\hat{s}}| + k_2\theta |h_s-y_{\wt{s}}| + k_2|h_s-x_s|\right)ds + M_t\\
&& + 3 k_3^2\int_{0}^{t} \frac{1}{m|h_s-x_s|}\left( (y_{\hat{s}})^{2q-1}y_s|sgn(z_s)-1|^2 + y_s|y_{\hat{s}}-y_{s}|^{2q-1} + |h_s-y_s| + |h_s-x_s|\right)ds \\
&\leq& k_2(1-\theta)\int_{0}^{t}|h_s-y_{\hat{s}}|ds + k_2\theta\int_{0}^{t} |h_s-y_{\wt{s}}|ds + \frac{3 k_3^2}{me_m}\int_{0}^{t} |h_s-y_s|ds +
 k_2\int_{0}^{t}|h_s-x_s|ds +M_t\\
&&+ \frac{3 k_3^2}{me_m}\int_{0}^{t} (y_{\hat{s}})^{2q-1}y_s|sgn(z_s)-1|^2ds + \frac{3 k_3^2}{me_m}\int_0^t y_s|y_{\hat{s}}-y_{s}|^{2q-1}ds
 +\frac{3 k_3^2T}{m},
\eeao
where in the second step we have used (\ref{eq23}) and (\ref{eq24}) and the properties of $\phi_m$ and
$$
M_t:= \int_{0}^{t} \phi_m^{\prime}(h_u-x_u)(sgn(z_u)g(y_{\hat{u}},y_{\wt{u}})-g(x_u,x_u)) dW_u.
$$
Taking expectations in the above inequality yields
\beao
&&\bfE\phi_m(h_{t}-x_{t})\leq k_2(1-\theta)\int_{0}^{t}\bfE|h_s-y_{\hat{s}}|ds + k_2\theta\int_{0}^{t}\bfE|h_s-y_{\wt{s}}|ds + \frac{3 k_3^2}{me_m}\int_{0}^{t}\bfE|h_s-y_s|ds\\
&&+\frac{3 k_3^2}{me_m}\int_{0}^{t}\bfE(y_{\hat{s}})^{2q-1}y_s|sgn(z_s)-1|^2ds + \frac{3 k_3^2}{me_m}\int_0^t \bfE y_s|y_{\hat{s}}-y_{s}|^{2q-1}ds
 +\frac{3 k_3^2T}{m} +  k_2\int_{0}^{t}\bfE|h_s-x_s|ds\\
&\leq&k_2(1-\theta)T\hat{C}_1\sqrt{\D} + k_2\theta T\wt{C}_1\sqrt{\D} + \frac{3 k_3^2TC_1}{me_m}\D +   k_2\int_{0}^{t}\bfE|h_s-x_s|ds\\
&&+\frac{3 k_3^2}{me_m}\int_{0}^{t}\bfE(y_{\hat{s}})^{2q-1}y_s|sgn(z_s)-1|^2ds + \frac{3 k_3^2}{me_m}\int_0^t \sqrt{\bfE(y_s)^2}\sqrt{\bfE|y_{\hat{s}}-y_{s}|^{4q-2}}ds  +\frac{3 k_3^2T}{m}\\
&\leq&k_2T((1-\theta)\hat{C}_1+\theta\wt{C}_1)\sqrt{\D} + \frac{3 k_3^2TC_1}{me_m}\D +   k_2\int_{0}^{t}\bfE|h_s-x_s|ds\\
&&+\frac{3 k_3^2}{me_m}\int_{0}^{t}\bfE(y_{\hat{s}})^{2q-1}y_s|sgn(z_s)-1|^2ds + \frac{3 k_3^2T}{me_m}\sqrt{A_2\bfE(x_0 + k_1T)^2}\sqrt{\hat{A}_{4q-2}}\D^{q-\frac{1}{2}}  +\frac{3 k_3^2T}{m},
\eeao
where we have used  Lemma \ref{l6} in the second step and Holder inequality, Lemmata \ref{l4} and \ref{l5} in the third step and the fact that $\bfE M_t=0$.\footnote{The function $d(u)=\phi_m^{\prime}(h_u-x_u)sgn(z_u)(g(y_{\hat{u}},y_{\wt{u}})-g(x_u,x_u))$ belongs to the space $\bbm^2([0,t];\bbR)$ of real valued measurable $\bbf_t-$adapted processes such that $\bfE\int_0^{t}|d(u)|^2du<\infty$ thus \cite[Th. 1.5.8]{mao:1997} implies $\bfE M_t=0$.} It holds for $s\in[t_n,t_{n+1}]$ that
\beao
&&\bfE(y_{\hat{s}})^{2q-1}y_s|sgn(z_s)-1|^2=\bfE\left(4(y_{t_n})^{2q-1}y_s\bbi_{\{z_s\leq0\}}\right)\\
\nonumber&\leq&4\bfE\left|(y_{t_n})^{2q-1}y_s-(y_{t_n})^{2q}\right|  + 4\bfE\left( (y_{t_n})^{2q}\bbi_{\{z_s\leq0\}}\bbi_{\{y_{t_n}\leq\D^{1-2\xi}\}}\right) + 4\bfE\left( (y_{t_n})^{2q}\bbi_{\{z_s\leq0\}}\bbi_{\{y_{t_n}>\D^{1-2\xi}\}}\right) \\
&\leq& 4\bfE\left|(y_{t_n})^{2q-1}(y_s - y_{t_n})\right|  + 4\D^{2q-4q\xi}  + 4\sqrt{\bfE (y_{t_n})^{4q}}\sqrt{\bfP\left({\{z_s\leq0\}}\cap{\{y_{t_n}>\D^{1-2\xi}\}}\right)}\\
&\leq& 4\sqrt{\bfE (y_{t_n})^{4q-2}}\sqrt{\bfE|y_s - y_{t_n}|^2} + 4\D^{2q-4q\xi}  + 4\sqrt{\bfE (y_{t_n})^{4q}}\sqrt{C_{k_2,k_3,\theta,\D}\sqrt{\D}}\\
&\leq&4\sqrt{A_{4q-2}\bfE(x_0+ k_1T)^{4q-2}}\sqrt{\hat{A}_2}\sqrt{\D} + 4\D^{2q-4q\xi}  + 4\sqrt{A_{4q}\bfE(x_0+ k_1T)^{4q}}\sqrt{C_{k_2,k_3,\theta,\D}}\D^{\frac{1}{4}},
\eeao
where we have used Holder inequality and Lemma \ref{l1}  in the third step and Lemmata \ref{l5} and \ref{l4} in the final step. For $\xi=\frac{1}{2}-\frac{1}{16q}$ we get the estimate
\beqq\label{eq25}
\bfE(y_{\hat{s}})^{2q-1}y_s|sgn(z_s)-1|^2\leq 4\left(\sqrt{A_{4q}\bfE(x_0+ k_1T)^{4q}C_{k_2,k_3,\theta,\D}}\vee1\right)\D^{\frac{1}{4}}.
\eeqq
Thus (\ref{eq22}) becomes
\beao
&&\bfE|h_{t}-x_{t}| \leq e_{m-1} + J_1\sqrt{\D} + 3 k_3^2TC_1\frac{\D}{me_m} + J_2\frac{\D^{\frac{1}{4}}}{me_m} + J_3\frac{\D^{q-\frac{1}{2}}}{me_m} +3 k_3^2T\frac{1}{m} + k_2\int_{0}^{t} \bfE|h_s-x_s|ds\\
&\leq& J_2\frac{\D^{\frac{1}{4}}}{me_m} + J_3\frac{\D^{q-\frac{1}{2}}}{me_m} +3 k_3^2T\frac{1}{m} + k_2\int_{0}^{t} \bfE|h_s-x_s|ds\\
&\leq&\left(J_2\frac{\D^{\frac{1}{4}}}{me_m} + J_3\frac{\D^{q-\frac{1}{2}}}{me_m} +3 k_3^2T\frac{1}{m}\right)e^{k_2t},
\eeao
where in the second step we have used the asymptotic relations, $ \D^{\ka}=o(\D^{q-\frac{1}{2}})$ as $\D\downarrow0$ for any $\ka\geq1/2, \,e_{m-1}=o(\frac{1}{m})$ as $\mto, \,\sqrt{\D}=o(\frac{\D^{\ka}}{me_m})$ for any $\ka\leq1$ as $\mto,$ in the last step we have used the Gronwall inequality and $J_3$ is as defined in Proposition \ref{pr1} while
$$J_1:=k_2T((1-\theta)\hat{C}_1+\theta\wt{C}_1).$$
Taking the supremum over all $0\leq t \leq T$ gives (\ref{eq21}).
\epf

\subsection{Convergence of the auxiliary process $(h_t)$ to $(x_t)$ in $\bbl^2$.}\label{ssec:s3.3}

\bpr\label{pr2}
Let Assumption A hold. Then we have
\beqq \label{eq27}
\bfE\sup_{0\leq t\leq T}|h_{t}-x_{t}|^2\leq\frac{C_{\ep}}{\sqrt{\ln (\D)^{-1}}},
\eeqq
where $C_{\ep}$ is independent of $\D$ and given by $C_{\ep}:=72\sqrt{\frac{3}{2\ep}}(k_3)^4T^2e^{6T^2k_2+k_2T},$ where $\ep$ is such that $0<\ep<\frac{1}{4}\wedge (q-\frac{1}{2}).$
\epr

\bpf[Proof of Proposition \ref{pr2}]
We estimate the difference $|\bbE_{t}|^2:=|h_{t}-x_{t}|^2.$ It holds that
\beao
&&|\bbE_{t}|^2=\left|\int_{0}^{t}\left(f(y_{\hat{s}},y_{\wt{s}})-f(x_s,x_s)\right)ds + \int_{0}^{t}\left(sgn(z_s)g(y_{\hat{s}},y_s)-g(x_s,x_s)\right)dW_s\right|^2\\
&\leq& 2T\int_{0}^{t}\left(k_2(1-\theta)|h_s - y_{\hat{s}}| + k_2\theta |h_s-y_{\wt{s}}| + k_2|\bbE_s|\right)^2ds + 2|M_t|^2\\
&\leq& 6Tk_2^2(1-\theta)^2\int_{0}^{t}|h_s - y_{\hat{s}}|^2ds + 6Tk_2^2\theta^2\int_{0}^{t}|h_s-y_{\wt{s}}|^2ds + 6Tk_2\int_{0}^{t}|\bbE_s|^2ds + 2|M_t|^2,
\eeao
where in the second step we have used Cauchy-Schwarz inequality and (\ref{eq23}) and
$$
M_t:= \int_{0}^{t} (sgn(z_u)g(y_{\hat{u}},y_{\wt{u}})-g(x_u,x_u)) dW_u.
$$
Taking the supremum over all $t\in[0,T]$ and then expectations we have
\beam
\nonumber
\bfE\sup_{0\leq t\leq T}|\bbE_t|^2&\leq& 6Tk_2^2(1-\theta)^2\int_{0}^{T}\bfE|h_s - y_{\hat{s}}|^2ds +
6Tk_2^2\theta^2\int_{0}^{T}\bfE|h_s-y_{\wt{s}}|^2ds\\
\nonumber&& + 6Tk_2\int_{0}^{T}\bfE\sup_{0\leq l\leq s}|\bbE_{l}|^2ds + 2\bfE\sup_{0\leq t\leq T} |M_t|^2\\
\label{eq26}&\leq& 6T^2k_2^2(1-\theta)^2\hat{A}_2\D + 6T^2k_2^2\theta^2\wt{A}_2\D + 6Tk_2\int_{0}^{T}\bfE\sup_{0\leq l\leq s}|\bbE_{l}|^2ds+ 8\bfE |M_T|^2,
\eeam
where in the second  step we have used Lemma \ref{l5} and Doob's martingale inequality with $p=2,$ since  $M_t$ is an $\bbR-$valued martingale that belongs to $\bbl^2.$ It holds that
\beao
&&\bfE |M_T|^2:=\bfE\left|\int_{0}^{T} |sgn(z_s)g(y_{\hat{s}},y_s)-g(x_s,x_s)|dW_s\right|^2\\
&=& \bfE\left(\int_{0}^{T} |sgn(z_s)g(y_{\hat{s}},y_s)-g(x_s,x_s)|^2ds\right)\\
&\leq& 3 k_3^2\bfE\left(\int_{0}^{T} \left( (y_{\hat{s}})^{2q-1}y_s(sgn(z_s)-1)^2 + y_s|y_{\hat{s}}-y_{s}|^{2q-1} + |h_s-y_s| + |h_s-x_s|\right)ds\right)\\
&\leq& 3k_3^2\int_{0}^{T} \bfE\left((y_{\hat{s}})^{2q-1}y_s|sgn(z_s)-1|^2\right)ds + 3k_3^2\int_{0}^{T}\bfE\left(y_s|y_{\hat{s}}-y_{s}|^{2q-1}\right)ds\\
&& + 3k_3^2\int_{0}^{T}\bfE|h_s-y_s|ds + 3k_3^2\int_{0}^{T}\bfE|\bbE_s|ds,
\eeao
where we have used (\ref{eq24}). Now, we use again the estimate (\ref{eq25}) to get
\beao
\bfE |M_T|^2&\leq& J_2\D^{\frac{1}{4}} + J_6\sqrt{\D^{2q-1}} + 3k_3^2TC_1\D + 3k_3^2\int_{0}^{T}\bfE|\bbE_s|ds\\
&\leq& J_2\D^{\frac{1}{4}} + J_6\D^{q-\frac{1}{2}} + 3k_3^2\int_{0}^{T}\bfE|\bbE_s|ds,
\eeao
where we have used the asymptotic relations, $\D^l=o(\D^{q-\frac{1}{2}})$ for all $l\geq\frac{1}{2}$ as $\D\downarrow0,$ the quantity $J_6$ is given by $J_6:=3k_3^2T\sqrt{A_2\bfE(x_0 + k_1T)^2}\sqrt{\hat{A}_{4q-2}}$ and $J_2$ is as  given in the statement of Proposition \ref{pr1} and depends on $\D$ through $C_{k_2,k_3,\theta,\D}$ where as already stated before we have that  $C_{k_2,k_3,\theta,\D}\rightarrow k_3,$ as $\D\downarrow0.$

Relation (\ref{eq26}) becomes
\beao
&&\bfE\sup_{0\leq t\leq T}|\bbE_t|^2\leq 8J_2\D^{\frac{1}{4}} +  8J_6\D^{q-\frac{1}{2}} +J_5\D + 6Tk_2\int_{0}^{T}\bfE\sup_{0\leq l\leq s}|\bbE_{l}|^2ds + 24k_3^2\int_{0}^{T}\bfE\sup_{0\leq l\leq s}|\bbE_{l}|ds\\
&\leq& 8J_2\D^{\frac{1}{4}} + 8J_6\D^{q-\frac{1}{2}} + 24k_3^2T\left(J_2\frac{\D^{\frac{1}{4}}}{me_m} + J_3\frac{\D^{q-\frac{1}{2}}}{me_m} +3 k_3^2T\frac{1}{m}\right)e^{k_2T} + 6Tk_2\int_{0}^{T}\bfE\sup_{0\leq l\leq s}|\bbE_{l}|^2ds\\
&\leq&  24k_3^2TJ_2e^{k_2T+6T^2k_2}\frac{\D^{\frac{1}{4}}}{me_m} +24k_3^2TJ_3e^{k_2T+6T^2k_2}\frac{\D^{q-\frac{1}{2}}}{me_m} + 72(k_3)^4T^2e^{k_2T+6T^2k_2}\frac{1}{m},
\eeao
where we have used Proposition \ref{pr1} in the second step with the sequence $e_m$ as defined there, Gronwall inequality in the last step and the asymptotic relation $\D^{\ka}=o(\frac{\D^{\ka}}{me_m})$ as $\mto,$ for any $\ka>0$ and $J_5$ is independent of $\D$ and given by
$J_5:=6T^2k_2^2((1-\theta)^2\hat{A}_2 +\theta^2\wt{A}_2).$

We take $m=\sqrt{\ln \D^{-\lam}},$ with $\lam>0$ to be specified soon and note that $e^{\sqrt{\ln \D^{-\lam}}}=o(\D^{-\lam})$, as $\D\downarrow0,$
since $e^{\sqrt{\ln n}}=o(n)$, as $\nto.$ Moreover it holds that
$$
\frac{\D^{q-\frac{1}{2}}}{e_m}=\frac{\D^{q-\frac{1}{2}}}{e^{-\frac{m^2}{2}}}e^{\frac{m}{2}}=\frac{\D^{q-\frac{1}{2}}}{e^{-\frac{\ln \D^{-\lam}}{2}}}e^{\frac{1}{2}\sqrt{\ln \D^{-\lam}}}=\D^{q-\frac{1}{2}-\frac{3\lam}{2}}\frac{e^{\frac{1}{2}\sqrt{\ln \D^{-\lam}}}}{\D^{-\lam}}.
$$
Now, since $q>\frac{1}{2}$ there is an $\ep>0$ small enough such that $q-\frac{1}{2}-\ep>0.$ We take $\lam=\frac{2\ep}{3}$ and conclude that
$$
\frac{\D^{q-\frac{1}{2}}}{e_m}=\D^{q-\frac{1}{2}-\ep}\frac{e^{\frac{1}{2}\sqrt{\ln \D^{-\frac{2\ep}{3}}}}}{\D^{-\frac{2\ep}{3}}}\rightarrow0,
$$
as $\D\rightarrow0$ which in turn implies the asymptotic relation $\frac{\D^{q-\frac{1}{2}}}{me_m}=o(\frac{1}{m})$ as $\mto,$ with the logarithmic rate stated before. In the same way we can show $\frac{\D^{\frac{1}{4}}}{me_m}=o(\frac{1}{m})$ as $\mto,$ by taking $\ep<\frac{1}{4}.$
We finally arrive at
$$
\bfE\sup_{0\leq t\leq T}|\bbE_t|^2\leq 72(k_3)^4T^2e^{k_2T+6T^2k_2}\frac{1}{\sqrt{\ln \D^{-\frac{2\ep}{3}}}},
$$
by taking $0<\ep<\frac{1}{4}\wedge (q-\frac{1}{2}),$ which implies (\ref{eq27}).
\epf

\subsection{Proof of Theorem \ref{t1}.}\label{ssec:s3.4}

In order to finish the proof of Theorem \ref{t1} we just use the triangle inequality, Lemma \ref{l6} and Proposition \ref{pr2} to get
\beao
\bfE\sup_{0\leq t\leq T}|y_{t}-x_{t}|^2 &\leq& 2\bfE\sup_{0\leq t\leq T}|h_{t}-y_{t}|^2 + 2\bfE\sup_{0\leq t\leq T}|\bbE_{t}|^2\\
&\leq& 2C_2\D^2 + 2\frac{C_{\ep}}{\sqrt{\ln \D^{-1}}}\leq \frac{C}{\sqrt{\ln \D^{-1}}},
\eeao
where $C=C(k_2,k_3,\ep,T),$ and given in the statement of Theorem \ref{t1}.

\section{Polynomial rate of convergence.}\label{sec:s4}
\setcounter{equation}{0}

We work with the stochastic time change inspired by \cite{berkaoui:2004}.
We define the process
$$
\gamma(t):=\int_0^t \frac{192 (k_3)^2q^2}{\left[(y_s)^{1-q} + (x_s)^{1-q}\right]^2}ds
$$
and the stopping time
$$
\tau_l:=\inf\{s\in[0,T]: 6T k_2s + \gamma(s)\geq l\}.
$$
The process $\gamma(t)$ is well defined since $x_t>0$ a.s. and $y_t\geq0$ (see Section \ref{sec:s2}).

The difference $|\bbE_{t}|^2:=|h_{t}-x_{t}|^2$ is estimated as in Section \ref{sec:s3} and we get,  as in (\ref{eq26}), that
\beqq\label{eq31}
\bfE\sup_{0\leq t\leq \tau}|\bbE_t|^2 \leq J_5\D + 6Tk_2\int_{0}^{\tau}\bfE\sup_{0\leq l\leq s}|\bbE_{l}|^2ds+ 8\bfE |M_\tau|^2,
\eeqq
where $\tau$ a stopping time and $J_5$ independent of $\D$ is as in proof of Proposition \ref{pr2}. The main difference here will be the estimation of the last term in (\ref{eq31}). The approach in Section \ref{sec:s3} resulted in the $\bbl^1$ estimation $\bfE|\bbE_{t}|$ where we used the Yamada-Watanabe approach.  Now, we use the Berkaoui approach. Relation (\ref{eq24}) becomes
\beao
&&|sgn(z_s)g(y_{\hat{s}},y_s)-g(x_s,x_s)|^2\leq 3 k_3^2 \left( (y_{\hat{s}})^{2q-1}y_s(sgn(z_s)-1)^2 + y_s|y_{\hat{s}}-y_{s}|^{2q-1} + |(y_s)^q-(x_s)^q|^2\right)\\
&\leq&3 k_3^2 \left( (y_{\hat{s}})^{2q-1}y_s(sgn(z_s)-1)^2 + y_s|y_{\hat{s}}-y_{s}|^{2q-1}\right) + |(y_s)^q-(x_s)^q|^2\left((y_s)^{1-q} + (x_s)^{1-q}\right)^2\frac{(\gamma_s)^{\prime}}{64q^2}\\
&\leq&3 k_3^2 \left( (y_{\hat{s}})^{2q-1}y_s(sgn(z_s)-1)^2 + y_s|y_{\hat{s}}-y_{s}|^{2q-1}\right) + \frac{1}{8}\left(|h_s - y_s|^2 + |\bbE_s|^2\right)(\gamma_s)^{\prime},
\eeao
where we have used the inequality
\beqq\label{eq32}
|a^q -b^q|(a^{1-q} + b^{1-q})\leq 2q|a-b|,
\eeqq
valid for all $a\geq0, b\geq0$ and $\frac{1}{2}\leq q\leq1.$ Consequently, we get the upper bound
\beao
&&\bfE |M_\tau|^2:=\bfE\left|\int_{0}^{\tau} |sgn(z_s)g(y_{\hat{s}},y_s)-g(x_s,x_s)|dW_s\right|^2\\
&\leq& 8J_2\D^{\frac{1}{4}} + J_6\D^{q-\frac{1}{2}} + \frac{1}{8}\int_{0}^{\tau}\bfE |h_s - y_s|^2(\gamma_s)^{\prime}ds +  \frac{1}{8}\int_{0}^{\tau}\bfE|\bbE_s|^2(\gamma_s)^{\prime}ds\\
&\leq& 8J_2\D^{\frac{1}{4}} + J_6\D^{q-\frac{1}{2}} + \frac{1}{8}\int_{0}^{\tau}\sqrt{\bfE |h_s - y_s|^4}\sqrt{\bfE((\gamma_s)^{\prime})^2}ds +  \frac{1}{8}\int_{0}^{\tau}\bfE|\bbE_s|^2(\gamma_s)^{\prime}ds,
\eeao
where we used the estimate (\ref{eq25}) and Holder inequality $J_2$ is as in the statement of Proposition \ref{pr1} and $J_6$ independent of $\D$ is as in proof of Proposition \ref{pr2}.
Relation (\ref{eq31}) becomes
\beao
&&\bfE\sup_{0\leq t\leq \tau}|\bbE_t|^2\leq 8J_2\D^{\frac{1}{4}} + 8J_6\D^{q-\frac{1}{2}} + 6Tk_2\int_{0}^{\tau}\bfE\sup_{0\leq l\leq s}|\bbE_{l}|^2ds + \int_{0}^{\tau}\sqrt{\bfE |h_s - y_s|^4}\sqrt{\bfE((\gamma_s)^{\prime})^2}ds\\
&&+ \int_{0}^{\tau}\bfE|\bbE_s|^2(\gamma_s)^{\prime}ds\\
&\leq& 8(J_2\vee J_6)\D^{(q-\frac{1}{2})\wedge\frac{1}{4}} + \sqrt{C_4}\D^2\int_{0}^{\tau}\sqrt{\bfE\left(\frac{192 (k_3)^2q^2}{\left[(y_s)^{1-q} + (x_s)^{1-q}\right]^2}\right)^2}ds +\int_{0}^{\tau}\bfE\sup_{0\leq l\leq s}|\bbE_l|^2(6Tk_2s + \gamma_s)^{\prime}ds\\
&\leq& 8(J_2\vee J_6)\D^{(q-\frac{1}{2})\wedge\frac{1}{4}} + \sqrt{C_4}192 (k_3)^2q^2\D^2\int_{0}^{\tau}\sqrt{\bfE\left(\frac{1}{(x_s)^{2-2q}}\right)}ds +\int_{0}^{\tau}\bfE\sup_{0\leq l\leq s}|\bbE_l|^2(6Tk_2s + \gamma_s)^{\prime}ds,
\eeao
where we have used Lemma \ref{l6} in the second step. At this point we want to estimate the inverse moments of $(x_t)$ and to do so we  consider the transformation $v=x^{2-2q}$ and apply Ito's formula to get
\beqq\label{eq45}
v_t=v_0 + \int_0^t \left( \un{(1-2q)(1-q)(k_3)^2}_{K_0} + \un{2(1-q)k_1}_{K_1}(v_s)^{\frac{1-2q}{2-2q}} - \un{2(1-q)k_2}_{K_2}v_s \right)ds + \int_{0}^{t}\un{2k_3(1-q)}_{K_3}\sqrt{v_s}dW_s,
\eeqq
for $t\in [0,T],$ where $v_0=(x_0)^{2-2q}>0.$ Denote the drift coefficient of the process $(v_t)$ by $a(v_t)$ and consider the function
\beqq\label{eq41}
\alpha(v):=a(v) - \lam + K_2v + \un{\frac{(2q-1)(\lam + K_0)^{\frac{1}{2q-1}}}{(k_1)^{\frac{2-2q}{2q-1}}}}_{\eta(\lam)}v,
\eeqq
where $\lam\geq0.$ Some elementary calculations show that this function attains its minimum at $v^*\!\!:=\!\!\left(\!\frac{k_1(2q-1)}{\eta(\lam)}\!\right)^{2-2q}$ and $\alpha(v^*)=0,$ thus
$$
a(v)\geq \lam - \left( K_2 + \eta(\lam)\right)v.
$$
Consider the process $(\zeta_t(\lam))$ defined through
\beqq\label{eq46}
\zeta_t(\lam)=\zeta_0 + \int_0^t (\lam - \left( K_2 + \eta(\lam)\right)\zeta_s)ds + \int_{0}^{t}K_3\sqrt{\zeta_s}dW_s,
\eeqq
for $t\in [0,T]$ with $\zeta_0(\lam)=v_0.$ Process (\ref{eq46}) is a square root diffusion process and when $\frac{2\lam}{(K_3)^2}-1\geq0$ or
\beqq\label{eq42}
\lam\geq 2(1-q)^2(k_3)^2,
\eeqq
is a CIR process which remains positive if $\zeta_0(\lam)>0.$ By a comparison theorem \cite[Prop. 5.2.18]{karatzas_shreve:1988} it holds that $v_t\geq \zeta_t(\lam)>0$ a.s. or $(v_t)^{-1}\leq(\zeta_t(\lam))^{-1}$ a.s. or equivalently $(x_t)^{2q-2}\leq(\zeta_t(\lam))^{-1}$ a.s.
The inverse moment bounds of $(\zeta_t(\lam))$ follow  by \cite[Rel. 3.1]{dereich_et_al:2011}
$$
\sup_{t\in[0,T]}\bfE (\zeta_t(\lam))^p<\infty, \qquad \mbox{ for }  p>-2\frac{\lam}{K_3^2}
$$
by choosing big enough $\lam$ and particularly such that (\ref{eq42}) holds strictly. Therefore,
\beqq\label{eq33}
\bfE\sup_{0\leq t\leq \tau}|\bbE_t|^2\leq 8(J_2\vee J_6)\D^{(q-\frac{1}{2})\wedge\frac{1}{4}} + \int_{0}^{\tau}\bfE\sup_{0\leq l\leq s}|\bbE_l|^2(6Tk_2s + \gamma_s)^{\prime}ds.
\eeqq
Relation (\ref{eq33}) for $\tau=\tau_l$ implies
\beam\nonumber
\bfE \sup_{0\leq t\leq \tau_l}(\bbE_t)^2 &\leq& 8(J_2\vee J_6)\D^{(q-\frac{1}{2})\wedge\frac{1}{4}}  + \int_0^{\tau_l} \bfE\sup_{0\leq l\leq s}(\bbE_l)^2 (6T k_2s + \gamma_s)^{\prime}ds\\
\nonumber&\leq& 8(J_2\vee J_6)\D^{(q-\frac{1}{2})\wedge\frac{1}{4}}  + \int_0^l \bfE\sup_{0\leq j\leq u}(\bbE_{\tau_j})^2du\\
\label{eq34} &\leq& 8(J_2\vee J_6)e^l\D^{(q-\frac{1}{2})\wedge\frac{1}{4}},
\eeam
where in the last step we have used Gronwall's inequality. Using again relation (\ref{eq33}) for $\tau=T$ and under the change of variables $u=6Tk_2s + \gamma_s$ we get,
\beao
&&\bfE \sup_{0\leq t\leq T}(\bbE_t)^2\leq 8(J_2\vee J_6)\D^{(q-\frac{1}{2})\wedge\frac{1}{4}} + \int_0^{6k_2T^2 + \gamma_T}\bfE\sup_{0\leq j\leq u}(\bbE_{\tau_j})^2du \\
&\leq&8(J_2\vee J_6)\D^{(q-\frac{1}{2})\wedge\frac{1}{4}} + \int_0^{\infty}\bfE\left(\sup_{0\leq j\leq u}(\bbi_{\{6k_2T^2 + \gamma_T\geq u\}}\bbE_{\tau_j})^2\right)du \\
&\leq&\int_0^{6k_2T^2}\bfE\sup_{0\leq j\leq u}(\bbE_{\tau_j})^2du + \int_{6k_2T^2}^{\infty}\bfP(6k_2T^2 + \gamma_T\geq u)\bfE\left(\sup_{0\leq j\leq u}(\bbE_{\tau_j})^2 \big| \{6k_2T^2 +\gamma_T\geq u\}\right)du\\
&&+ 8(J_2\vee J_6)\D^{(q-\frac{1}{2})\wedge\frac{1}{4}}\\
&\leq&8(J_2\vee J_6)e^{6k_2T^2}\D^{(q-\frac{1}{2})\wedge\frac{1}{4}} + \int_0^{\infty}\bfP(\gamma_T\geq u)\bfE\sup_{0\leq j\leq u}(\bbE_{\tau_j})^2 du + 8(J_2\vee J_6)\D^{(q-\frac{1}{2})\wedge\frac{1}{4}}\\
&\leq&16(J_2\vee J_6)e^{6k_2T^2}\D^{(q-\frac{1}{2})\wedge\frac{1}{4}} + 8(J_2\vee J_6)\D^{(q-\frac{1}{2})\wedge\frac{1}{4}}\int_0^{\infty}\bfP(\gamma_T\geq u)e^udu,
\eeao
where in the last steps we have used (\ref{eq34}). We proceed by showing that $u\mapsto \bfP(\gamma_T\geq u)e^u\in\bbl^1(\bbR_+).$ It holds that
\beqq \label{eq35}
 \bfP(\gamma_T\geq u)\leq e^{-\ep u}\bfE(e^{\ep \gamma_T}),
\eeqq
 for any $\ep>0$ by Markov inequality. The following bound holds
$$
\gamma_T=\int_0^T \frac{192 (k_3)^2q^2}{\left[(y_s)^{1-q} + (x_s)^{1-q}\right]^2}ds\leq 192 (k_3)^2q^2 \int_0^T (x_s)^{2q-2}ds,
$$
thus
\beqq \label{eq44}
\bfE(e^{\ep\gamma_T})\leq \bfE\left(e^{\ep 192 (k_3)^2q^2\int_0^T (x_s)^{2q-2}ds} \right),
\eeqq
where $-1<2q-2<0.$ It remains to bound the exponential inverse moments of $(x_t)$ defined through the stochastic integral equation (\ref{eq2}). Exponential inverse moments for the CIR process are known  \cite[Th. 3.1]{hurd_kuznetsov:2008} and are given by
\beqq\label{eq47}
\bfE e^{\del\int_0^t(\zeta_s(\lam))^{-1}ds}\leq C_{HK}(\zeta_0)^{-\frac{1}{2}(\nu(\lam)-\sqrt{\nu(\lam)^2+8 \frac{\del}{(K_3)^2}})},
\eeqq
for $0\leq \del\leq \left(\frac{2\lam}{K_3^2} - 1\right)^2\frac{K_3^2}{8}=:\nu(\lam)^2\frac{K_3^2}{8},$ where the positive constant $C_{HK}$ is explicitly given in \cite[Rel. (10)]{hurd_kuznetsov:2008} depends on the parameters $k_2,k_3,T,q,$ but is independent of $\zeta_0.$ Thus the other condition that we require for parameter $\lam$ is
\beqq\label{eq43}
\lam\geq 2(1-q)\sqrt{2\del}(k_3) + 2(1-q)^2(k_3)^2.
\eeqq
When (\ref{eq43}) is satisfied then (\ref{eq42}) is satisfied too, thus there is actually no restriction on the coefficient $\del$ in (\ref{eq47}) since we can always choose appropriately a $\lam$ such that  (\ref{eq43})  holds. Relation (\ref{eq44}) becomes
\beqq\label{eq48}
\bfE(e^{\ep\gamma_T})\leq \bfE\left(e^{\ep 192 (k_3)^2q^2\int_0^T (v_s)^{-1}ds} \right)\leq \bfE\left(e^{\ep 192 (k_3)^2q^2\int_0^T (\zeta_s(\lam))^{-1}ds} \right).
\eeqq
We therefore require that
\beqq\label{eq48.1}
192 (k_3)^2q^2\ep\leq \left(\nu(\lam)\right)^2\frac{K_3^2}{8}
\eeqq
and can always find a $\ep>1,$ such the above relation holds by choosing appropriately $\lam$ as discussed before.
Relation (\ref{eq48}) becomes
$$
\bfE(e^{\ep\gamma_T})\leq C_{HK}(\zeta_0)^{-\frac{\nu(\lam)}{2}},
$$
and therefore
$$
\bfP(\gamma_T\geq u)\leq C_{HK}(x_0)^{(1-q)\nu(\lam)}e^{-\ep u},
$$
where $\lam$ is chosen such that (\ref{eq48.1}) holds with $\ep>1.$
We conclude
\beao
\bfE \sup_{0\leq t\leq T}(\bbE_t)^2 &\leq& 16(J_2\vee J_6)e^{6k_2T^2}\D^{(q-\frac{1}{2})\wedge\frac{1}{4}}   + 8(J_2\vee J_6)C_{HK}(x_0)^{(1-q)\nu(\lam)}\D^{(q-\frac{1}{2})\wedge\frac{1}{4}}\int_0^{\infty}e^{(1-\ep) u}du\\
&\leq& C\cdot\D^{(q-\frac{1}{2})\wedge\frac{1}{4}},
\eeao
by choosing  $\ep>1,$ where $C=C(k_1,k_2,k_3,T,q,\ep):=8(J_2\vee J_6)\left(2e^{6k_2T^2} + \frac{C_{HK}}{\ep-1}(x_0)^{(1-q)\nu(\lam)}\right),$  is as given in  statement of Theorem \ref{t2}.

\section{Alternative approach improving the rate of convergence.}\label{sec:s5}
\setcounter{equation}{0}

\bpf[Proof of Theorem \ref{t3}]
We will discuss the proof and highlight the differences since one can follow the proofs of Theorems \ref{t1} and \ref{t2}. First of all note that Lemmata \ref{l4}, \ref{l5} and \ref{l6} still hold, i.e. the moment bounds and error bounds of $(y_t^{SD}),$ as well as the moment bounds involving the auxiliary process $(h_t)$ are true. The error $\bbE_t:=h_t-x_t$ now reads as
\beqq\label{eq106}
h_t-x_t=\int_0^t \left(f_{\theta}(y_{\hat{s}},y_{\wt{s}})-f_{\theta}(x_s,x_s)\right)ds + \int_{0}^{t}\left(g(y_{\hat{s}},y_s)-g(x_s,x_s)\right)dW_s.
\eeqq

As regards the first part of Theorem \ref{t3}, we now have that
\beqq \label{eq107}
\sup_{0\leq t\leq T}\bfE|\bbE_{t}|\leq\left(J_3\frac{\D^{q-\frac{1}{2}}}{me_m} +2(k_3)^2T\frac{1}{m}\right)e^{k_2T},
\eeqq
for any $m>1,$ where $e_m=e^{-m(m+1)/2}$ and $J_3$ as stated in Proposition \ref{pr1}. The bound (\ref{eq107}) follows in the same lines as in the proof of Proposition \ref{pr1} where now (\ref{eq24}) becomes
\beqq\label{eq108}
|g(y_{\hat{s}},y_s)-g(x_s,x_s)|^2\leq 2(k_3)^2 \left( y_s|y_{\hat{s}}-y_{s}|^{2q-1} + |h_s-y_s| + |h_s-x_s|\right)
\eeqq
and the term (\ref{eq25}) disappears. Moreover, we have that
\beqq \label{eq109}
\bfE\sup_{0\leq t\leq T}|\bbE_{t}|^2\leq\frac{C_{\ep}}{\sqrt{\ln (\D)^{-1}}},
\eeqq
where $C_{\ep}$ is independent of $\D$ and given by $C_{\ep}:=32\sqrt{\frac{3}{2\ep}}(k_3)^4T^2e^{6T^2k_2+k_2T},$ where $\ep$ is such that $0<\ep<q-\frac{1}{2}.$
The bound (\ref{eq109}) follows in the same lines as in the proof of Proposition \ref{pr2} where now $M_t:= \int_{0}^{t} (g(y_{\hat{u}},y_{\wt{u}})-g(x_u,x_u)) dW_u$ and
$$\bfE |M_T|^2\leq 2(k_3)^2\int_{0}^{T}\bfE\left(y_s|y_{\hat{s}}-y_{s}|^{2q-1}\right)ds + 2(k_3)^2\int_{0}^{T}\bfE|h_s-y_s|ds + 2(k_3)^2\int_{0}^{T}\bfE|\bbE_s|ds,
$$
implying
$$\bfE\sup_{0\leq t\leq T}|\bbE_t|^2\leq \frac{32}{3}(k_3)^2TJ_3e^{k_2T+6T^2k_2}\frac{\D^{q-\frac{1}{2}}}{me_m} + 32(k_3)^4T^2e^{k_2T+6T^2k_2}\frac{1}{m},$$
which in turn gives
$$
\bfE\sup_{0\leq t\leq T}|\bbE_t|^2\leq 32(k_3)^4T^2e^{k_2T+6T^2k_2}\frac{1}{\sqrt{\ln \D^{-\frac{2\ep}{3}}}},
$$
where now we take  $0<\ep<q-\frac{1}{2}.$

As regards the second part of Theorem \ref{t3}, we follow Section \ref{sec:s4}, where now we use the process
$$
\gamma(t):=\int_0^t \frac{128(k_3)^2q^2}{\left[(y_s)^{1-q} + (x_s)^{1-q}\right]^2}ds
$$
and the estimate
$$|g(y_{\hat{s}},y_s)-g(x_s,x_s)|^2\leq 2 (k_3)^2y_s|y_{\hat{s}}-y_{s}|^{2q-1} + \frac{1}{8}\left(|h_s - y_s|^2 + |\bbE_s|^2\right)(\gamma_s)^{\prime},
$$
to get the upper bound
$$
\bfE |M_\tau|^2\leq\frac{2}{3}J_6\D^{q-\frac{1}{2}} + \frac{1}{8}\int_{0}^{\tau}\sqrt{\bfE |h_s - y_s|^4}\sqrt{\bfE((\gamma_s)^{\prime})^2}ds +  \frac{1}{8}\int_{0}^{\tau}\bfE|\bbE_s|^2(\gamma_s)^{\prime}ds,
$$
which in turn implies first
$$\bfE\sup_{0\leq t\leq \tau_l}(\bbE_t)^2\leq \frac{16}{3}J_6e^l\D^{(q-\frac{1}{2})}$$
and then
$$\bfE \sup_{0\leq t\leq T}(\bbE_t)^2\leq \frac{32}{3}J_6e^{6k_2T^2}\D^{(q-\frac{1}{2})} + \frac{16}{3}J_6\D^{(q-\frac{1}{2})}\int_0^{\infty}\bfP(\gamma_T\geq u)e^udu.
$$
Finally, in order to bound $\bfE\left(e^{\ep 128 (k_3)^2q^2\int_0^T (x_s)^{2q-2}ds} \right),$ we require that
\beqq\label{eq110}
128 (k_3)^2q^2\ep\leq \left(\nu(\lam)\right)^2\frac{K_3^2}{8}
\eeqq
and can always find a $\ep>1,$ such that (\ref{eq110}) holds yielding
$$\bfE \sup_{0\leq t\leq T}(\bbE_t)^2 \leq C\cdot\D^{(q-\frac{1}{2})},
$$
where $C=C(k_1,k_2,k_3,T,q,\ep):=\frac{16}{3}J_6\left(2e^{6k_2T^2} + \frac{C_{HK}}{\ep-1}(x_0)^{(1-q)\nu(\lam)}\right),$ is as given in  statement of Theorem \ref{t3}.
\epf

\section{Numerical Experiments.}\label{sec:s6}
\setcounter{equation}{0}

We discretize the interval $[0,T]$ with a number of steps in power of $2.$ The semi-discrete (SD) scheme is given by
\beqq\label{eq74}
y_{t_{n+1}}^{SD}=\left( \sqrt{y_{t_n}\left(1- \frac{k_2\D}{1+k_2\theta\D}\right) + \frac{k_1\D}{1+k_2\theta\D} - \frac{(k_3)^2\D}{4(1+k_2\theta\D)^2}(y_{t_n})^{2q-1}} +\frac{k_3}{2(1+k_2\theta\D)}(y_{t_n})^{q -\frac{1}{2}}\D W_n\right)^2,
\eeqq
for $n=0,\ldots,N-1,$ where $\D W_n:=W_{t_{n+1}}-W_{t_{n}}$ are the increments of the Brownian motion which are standard normal r.v's.

The ALF (Alfonsi) scheme \cite[Sec. 3]{alfonsi:2013} is an implicit scheme which requires solving the nonlinear equation
\beqq\label{eq73}
Y_{n+1}= y_{t_n} + (1-q)\left(k_1(Y_{n+1})^{\frac{-q}{1-q}} - k_2Y_{n+1}- \frac{q(k_3)^2}{2}(Y_{n+1})^{-1}\right)\D + k_3(1-q)\D W_n,
\eeqq
and then computing $y_{t_{n+1}}^{ALF}=(Y_{n+1})^{\frac{1}{1-q}}.$ The estimation of $Y_{n+1}$ in (\ref{eq73}) can be done for example with Newton's method, but requires a small enough $\D.$\footnote{In the CIR case, i.e. when $q=1/2$ (\ref{eq73}) simplifies to a solution of a quadratic equation.}
We also consider a scheme recently proposed in
\cite{halidias:2015b} using again the SD method, but in a
different way,
\beqq\label{eq75} y_{t_{n+1}}^{HAL}(q)= \Big|
\left(y_{t_n}(1- k_2\D) + k_1\D  -
\frac{q(k_3)^2\D}{2}(y_{t_n})^{2q-1}\right)^{1-q} + k_3(1-q)\D
W_n\Big|^{\frac{1}{1-q}}, \eeqq 
for $n=0,\ldots,N-1.$ Note the similarity in the expressions of (\ref{eq75}) and the SD scheme (\ref{eq74}) proposed here. This  is not strange, because they
both rely in the same way of splitting the drift coefficient. In
particular, in the explicit HAL scheme, the following process is
considered
$$
y_t^{HAL}(q)=y_{t_n} + \wt{f}_1(y_{t_n})\cdot\D + \int_{t_n}^t\wt{f}_2(y_{s})ds + \int_{t_n}^{t} sgn(z_s)\wt{g}(y_s)dW_s,
$$
for $t\in(t_n,t_{n+1}]$ with $y_0=x_0$ a.s. where now
\beqq\label{eq76}
f(x)= \underbrace{k_1 - k_2x - \frac{q(k_3)^2}{2}x^{2q-1}}_{\wt{f}_1(x)} + \underbrace{\frac{q(k_3)^2}{2}x^{2q-1}}_{\wt{f}_2(x)}. \quad \wt{g}(x)=k_3 x^{q}
\eeqq
and
\beqq\label{eq77}
z_t=\left(y_{t_n}(1- k_2\D) + k_1\D  - \frac{q(k_3)^2\D}{2}(y_{t_n})^{2q-1}\right)^{1-q} + k_3(1-q)(W_t-W_{t_n}).
\eeqq
A comparison  with (\ref{eq28}) and (\ref{eq29}) shows $\wt{f}_2(x)=2qf_2(x)$ and $\wt{g}(x)=g(x,x),$ for $\theta=0.$ We write (\ref{eq76}) again as
\beqq\label{eq78}
y_t^{HAL}(q)=y_{t_n} + \left(k_1 - k_2y_{t_n} - \frac{q(k_3)^2}{2}(y_{t_n})^{2q-1}\right)\cdot\D
 + \int_{t_n}^{t}\frac{q(k_3)^2}{2}(y_{s})^{2q-1} ds + k_3\int_{t_n}^{t} sgn(z_s)(y_s)^qdW_s
\eeqq
and the process (\ref{eq78}) is well defined when
\beqq\label{eq79}
(k_3)^2\leq \frac{2}{q}k_1 \, \mbox{ and }  \D\leq\frac{2}{2k_2+q(k_3)^2}.
\eeqq
The reader can compare again with (\ref{eq3}) for $\theta=0.$ Solving for $y_t,$ we end up with $y_t^{HAL}(q)=|z_t|^{\frac{1}{1-q}}.$
The main result in \cite{halidias:2015b} is
$$
\bfE |y_t^{HAL} - x_t|^2\leq C\cdot \D^{2q(q-\frac{1}{2})},
$$
when (\ref{eq79}) holds, implying a rate of convergence at least $q(q-\frac{1}{2})$ which is bigger than the rate of convergence of the SD scheme proposed here which is at least $\frac{1}{2}(q-\frac{1}{2})$ (see Th.\ref{t3}).

We also consider two more linear-implicit schemes that were stated in the introduction and discussed in Appendix \ref{A1}. Namely, we compare with the balanced implicit method (BIM) with appropriate weight functions to guarantee positivity (\cite[Th. 5.9]{kahl_schurz:2006}), which reads
\beqq\label{eq80} 
y_{t_{n+1}}^{BIM}(q)= \frac{ y_{t_n} + k_1 \D +  k_3(y_{t_n})^{q}(\D W_n + |\D W_n|)}{ 1 + k_2\D + k_3 (y_{t_n})^{q-1}|\D W_n|}, \eeqq
and the balanced Milstein method (BMM) with the suggested weight functions \cite[Th. 5.9]{kahl_schurz:2006} that is given by
\beqq\label{eq81} 
y_{t_{n+1}}^{BMM}(q)= \frac{ y_{t_n} + (k_1 + (\Theta-1)k_2y_{t_n})\D +  k_3(y_{t_n})^{q}\D W_n +  \frac{q(k_3)^2}{2}(y_{t_n})^{2q-1}(\D W_n)^2}{ 1 + \Theta k_2\D + 
\frac{q(k_3)^2}{2}|y_{t_n}|^{2q-2}\D}. \eeqq
We take the relaxation parameter $\Theta$ to be $1/2$ as recommended  in \cite[Rel. 5.10]{kahl_schurz:2006}.

We aim to show experimentally the order of convergence for the above  positivity preserving methods for the estimation of the true solution of the CEV model (\ref{eq2}), i.e the semi-discrete methods SD method (\ref{eq74}) and the HAL scheme  (\ref{eq75}), as well as the implicit ALF scheme (\ref{eq73}) and the linear-implicit schemes BIM and BMM. The choice of the parameters is the same  as in \cite[Fig. 6]{kahl_jackel:2006} with $k_3=0.4.$
In particular $(x_0, k_1, k_2, k_3, q, T)=(\frac{1}{16}, \frac{1}{16}, 1, 0.4, \frac{3}{4}, 1).$

Furthermore, we would also like to reveal the dependence of the order of the semi discrete methods on $q,$ i.e. we want to verify our theoretical results and in particular the order shown in Theorem \ref{t3}. We take the level of implicitness of SD method (\ref{eq74})  to be $\theta=1,$ i.e. we consider the fully implicit scheme. We also discuss about the fully  explicit scheme, i.e. when $\theta=0,$ but also an intermediate scheme $\theta=1/2,$ in Section \ref{sec:s7}.

We want to estimate the endpoint $\bbl^2-$norm $\ep=\sqrt{\bfE|y^{(\D)}(T) - x_T|^2},$ of the difference between the numerical scheme evaluated at step size $\D$ and the exact solution of (\ref{eq2}). For that purpose, we 
compute $M$ batches of $L$ simulation paths, where each batch is estimated by 
$\hat{\ep_j}=\frac{1}{L}\sum_{i=1}^L|y_{i,j}^{(\D)}(T) - y_{i,j}^{(ref)}(T)|^2$ 
and the Monte Carlo estimator of the error is
\beqq\label{l2error}
\hat{\ep}=\sqrt{\frac{1}{ML}\sum_{j=1}^M\sum_{i=1}^L|y_{i,j}^{(\D)}(T) - y_{i,j}^{(ref)}(T)|^2},
\eeqq
and requires $M\cdot L$ Monte Carlo sample paths. The reference solution is evaluated at step size $2^{-14}$  of the numerical scheme. For the SD case, we have shown in Theorems \ref{t1}, \ref{t2} and \ref{t3}, that it strongly converges to the exact solution.
We simulate $100\cdot 100=10000$ paths, where the choice for $L=100$ is as  in \cite[p.118]{kloeden_platen_schurz:2003}. The choice of the number of trajectories  $M\cdot L=10^4$ is also considered in \cite[Sec.5]{tretyakov_zhang:2013} where a fundamental mean-square theorem is proved for SDEs with superlinear growing coefficients satisfying a one-side Lipschitz condition, but unfortunately it is not positivity preserving. Of course, the number of Monte Carlo paths has to be sufficiently large, so as not to significantly hinder the mean square errors. 

\begin{figure}[ht]
  \caption{ \small Convergence of fully implicit SD, HAL, ALF, BIM and BMM schemes applied to SDE (\ref{eq2}) with parameters $(x_0, k_1, k_2, k_3, q, T)=(\frac{1}{16}, \frac{1}{16}, 1, 0.4, \frac{3}{4}, 1)$ with $32$ digits of accuracy.}
  \centering
   \includegraphics[width=0.5\textwidth]{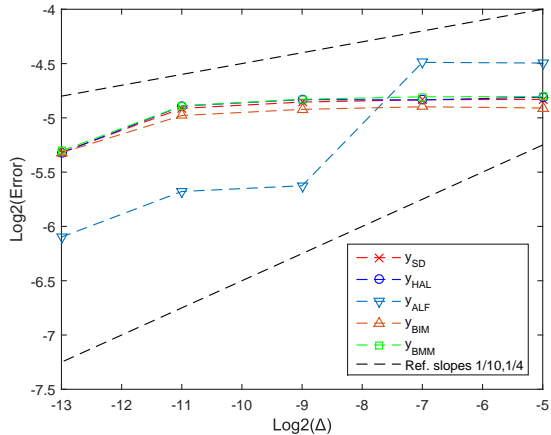}\label{SD_HAL_ALF_BIM_BMM_0.75}
  \end{figure}

We plot in a $\log_2-\log_2$ scale and error bars represent  $98\%-$confidence intervals. The results are shown in  Table \ref{tab:errors_0.75} and Figure \ref{SD_HAL_ALF_BIM_BMM_0.75}.
 Table \ref{tab:errors_0.75} does  not present the computed Monte Carlo errors with $98\%$ confidence, since they were at least $9$ times smaller that the mean-square errors.  

\begin{table}[htbp]  \scriptsize
\centering
        \begin{tabular}{|c|l|l|l|l|l|}
        \hline  Step $\D$ & $98\%$ SD-Error$(\theta=1)$ &  $98\%$ HAL-Error & $98\%$ ALF-Error & $98\%$ BIM-Error& $98\%$ BMM-Error\\
 \hline $2^{-5}$  & $0.0351607273610989$ &$0.0356800652730259$    &$0.0443426636440914$  & $0.0332754182024868$  & $0.0358194304596254$ \\
\hline $2^{-7}$  &  $0.0350820054172969$  & $0.035033855267215$   & $0.0445388652736149$ & $0.0335414964013289$  & $0.035736599920736$ \\
\hline $2^{-9}$  &  $0.0345654286067145$  &  $0.0351090011639902$  & $0.0202607612841043$ & $0.0329881367906935$  & $0.0351924591356631$ \\
\hline $2^{-11}$  & $0.0332045173200957$   & $0.0337092293915926$   & $0.0195315499600022$ & $0.0317244587795127$  & $0.0337311936077772$ \\
\hline $2^{-13}$  & $0.0250782316352445$  &  $0.0249779540239864$  & $0.0146002661407415$ & $0.0249983526384181$  & $0.025291682868944$ \\
 \hline  \end{tabular}
    \caption{\small Error and step size of fully implicit SD, HAL,ALF, BIM and BMM scheme for (\ref{eq2}) with $(x_0, k_1, k_2, k_3, q, T)=(\frac{1}{16}, \frac{1}{16}, 1, 0.4, \frac{3}{4}, 1)$  and  $32$ digits of accuracy.}
    \label{tab:errors_0.75}
\end{table}

\begin{figure}[ht]
  \caption{ \small Strong convergence error of the mean-reverting CEV process (\ref{eq2}) as a function of CPU time (in sec) using positivity preserving schemes SD,HAL, ALF, BIM and BMM with $(x_0,k_1, k_2, k_3, q, T)=(\frac{1}{16}, 100, 0.05, \frac{1}{16}, 1, 0.4, \frac{3}{4}, 1),$ and  $32$ digits of accuracy.}
  \centering
   \includegraphics[width=0.5\textwidth]{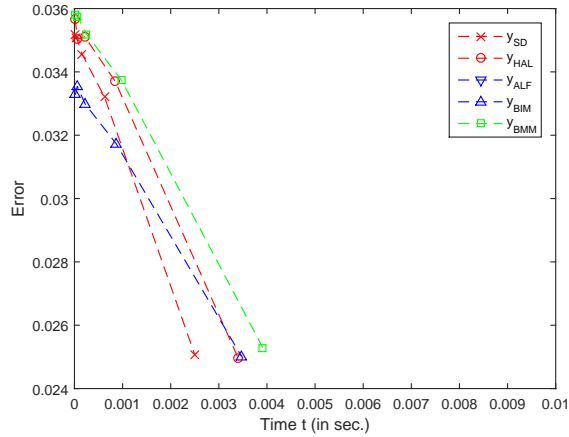}\label{CEVModel_075wrtTIME}
  \end{figure}

In Table \ref{tab:order_time_SD_0.75} we present the computational times,\footnote{\label{simu}We simulate with $3.06$GHz Intel Pentium, $1.49$GB of RAM in Matlab $R2014b$ Software. The random number generator is Mersenne Twister. The evaluated times do not include the random number generation time, since all the methods we compare, involve the same amount of random numbers.}  of fully implicit SD, HAL, ALF, BIM and BMM, for the same problem. Figure \ref{CEVModel_075wrtTIME}
shows the relation between the error and computer time consumption. As one can see from Table \ref{tab:order_time_SD_0.75} the CPU times for ALF are at least $1000$ times bigger than the other schemes, thus we choose in  Figure \ref{CEVModel_075wrtTIME} to restrict our attention to the rest of the methods.
	
		\begin{table}[htbp] \small 
\centering
        \begin{tabular}{|c|c|c|c|c|c|c|}
        \hline          $q=0.75$ & Step $\D$ & Implicit SD  & HAL& ALF & BIM& BMM\\
        \hline Time/Path(in sec) &$2^{-5}$ & $0.000013$  &$0.0000164$ & $0.0221883$ &$0.0000174$ &$0.0000196$ \\
        \hline Time/Path(in sec) &$2^{-7}$ & $0.0000422$ &$0.0000558$ & $0.0841705$ &$0.0000584$ &$0.0000657$ \\
         \hline Time/Path(in sec) &$2^{-9}$& $0.0001586$  &$0.0002137$& $0.2453943$ &$0.0002207$ &$0.0002482$ \\
        \hline Time/Path(in sec) &$2^{-11}$& $0.0006243$  &$0.0008437$& $0.9768619$ &$0.0008703$ &$0.0009795$\\      
        \hline Time/Path(in sec) &$2^{-13}$& $0.0024975$  &$0.0033977$& $3.9096332$ &$0.0034785$ &$0.0039143$ \\   
    \hline
  \end{tabular}
    \caption{\small Average computational time (in seconds) for a path, for different discretizations,  for all considered positivity preserving methods for the mean-reverting CEV process (\ref{eq2}) with $q=0.75.$}
    \label{tab:order_time_SD_0.75}
\end{table}

Finally, Table \ref{tab:orderall} presents the exact values of order of convergence for SD, HAL,ALF, BIM and BMM  produced by linear regression with the method of least squares fit, in the case one considers in one case $3$ points with steps $\D=2^{-9}, 2^{-11}, 2^{-13}$ and in other case all $5$ points including $\D=2^{-5}, 2^{-7}.$
    \begin{table}[htbp] \small
\centering
        \begin{tabular}{|c|c|c|c|c|c|c|}
        \hline  $q$ & $N^o$ of Points & SD$(\theta=1)$ & HAL & ALF& BIM & BMM$(\Theta=\frac{1}{2})$ \\
    \hline   $0.75$    & $5$ & $0.053$  & $0.054$ & $0.220$  & $0.045$& $0.054$\\
        \cline{2-7}   &  $3$ & $0.116$  & $0.123$ & $0.118$  & $0.100$ & $0.119$\\
        \hline
  \end{tabular}
    \caption{\small Order of convergence of SD with $\theta=1$, HAL, ALF, BIM and BMM approximation of (\ref{eq2}) with $(x_0, k_1, k_2, k_3, q, T)=(\frac{1}{16}, \frac{1}{16}, 1, 0.4, \frac{3}{4}, 1).$}
    \label{tab:orderall}
\end{table}

We show below, in Table \ref{tab:dif_CEV_0.75}, the $\bbl^2-$distance between our proposed method and the other methods for the numerical approximation of (\ref{eq2}). We work as before and estimate  the distance
\beqq\label{l2dist}
d(G,H)=\sqrt{\frac{1}{ML}\sum_{j=1}^M\sum_{i=1}^L|y_{i,j}^{(\D,G)}(T) - y_{i,j}^{(\D,H)}(T)|^2},
\eeqq
between method $G$ and $H,$ by considering  sufficient small $\D,$ and in particular for $\D=10^{-2}, 10^{-3}, 10^{-4}.$ 

    \begin{table}[htbp] \scriptsize
\centering
        \begin{tabular}{|c|l|l|l|l|}
        \hline  Step $\D$ & $98\%-d(SD,HAL)$  &  $98\%-d(SD,ALF)$       &  $98\%-d(SD,BIM)$ & $98\%-d(SD,BMM)$\\
        \hline      $10^{-2}$ & $0.000572742828312345$    & $0.0716139598501867$   & $0.0038372799928164$ & $0.000531186647260291$\\
        \hline      $10^{-3}$ & $0.000157661207152358$    & $0.0286630459661306$   & $0.00134601166557609$& $0.000156449754950229$\\
        \hline      $10^{-4}$ & $0.000049813956973071$    & $0.0283116960890337$   & $0.00044483884745991$& $0.000049771403122846$\\
 \hline
  \end{tabular}
    \caption{\small The $\bbl^2-$distance  between all the considered numerical schemes applied to SDE (\ref{eq2}) with parameter set $(x_0, k_1, k_2, k_3, q, T)=(\frac{1}{16}, \frac{1}{16}, 1, 0.4, \frac{3}{4}, 1)$.}
    \label{tab:dif_CEV_0.75}
\end{table}


Finally, we examine the behavior of SD w.r.t the parameter $q$. We want to examine the impact of $q$ in the order of convergence and verify our theoretical results and in particular Theorem \ref{t3}. Table \ref{tab:orderSD} shows the order of SD w.r.t $q.$

    \begin{table}[htbp] \small
\centering
        \begin{tabular}{|c|c|}
        \hline  $q$ &  Order of Fully Implicit SD\\
       \hline   $0.6$ & $0.110\,\,(0.052)$\\
        \hline  $0.7$ & $0.118\,\,(0.0525)$\\
        \hline   $0.8$ & $0.116\,\,(0.05)$  \\
        \hline   $0.9$ & $0.121\,\,(0.053)$ \\
        \hline
  \end{tabular}
    \caption{\small Order of convergence of SD$(\theta=1)$ approximation of (\ref{eq2}) with $(x_0, k_1, k_2, k_3, T)=(\frac{1}{16}, \frac{1}{16}, 1, 0.4, 1)$ for different values of $q$.}
    \label{tab:orderSD}
\end{table}

The following points of discussion are worth mentioning.
\begin{itemize}
    \item The performance of all methods, as shown in Table \ref{tab:errors_0.75} and  Figure \ref{SD_HAL_ALF_BIM_BMM_0.75}, implies, in terms of error estimates, that the  implicit ALF scheme performs better, for values of discretization steps $\D\geq 2^{-9}.$ All the other methods, i.e. the semi discrete SD and HAL, and the BIM and BMM have a similar behavior for all $\D$'s in the sense above as Figure \ref{SD_HAL_ALF_BIM_BMM_0.75} shows. The similarity of SD, HAL, BIM and BMM is also indicated in Table \ref{tab:dif_CEV_0.75}, where we see how close they are w.r.t. the $\bbl^2-$norm,
 and in Table \ref{tab:orderall}	where the convergence order is considered. Nevertheless,  Table \ref{tab:dif_CEV_0.75} also shows that in order to get an accuracy to at least $2$ decimal digits, which in practice may be adequate concerning that we want for example to evaluate an option and thus our results are in euros, there is no actual harm in choosing whatever of the above available methods. We may then choose the fastest one, as will be discussed later on.
    \item  We see that the strong order of convergence of implicit SD for  problem (\ref{eq2}) is at least $1/2(q-1/2)=1/8,$ as shown theoretically and presented in Table \ref{tab:orderall}. We also see that all methods converge with similar orders and the theoretically rate $1$ of the ALF method \cite{alfonsi:2013} does not hold for these $\D$'s. Thus, again we see that the rate in practical situations does not necessarily matter, if one has to compute with very small  $\D$'s to achieve it. Moreover, we present in Table \ref{tab:exp_SD_0.75}   the performance of the explicit SD method and see that it is very close to the implicit, which is of course natural to happen.
    \begin{table}[htbp] \small
\centering
        \begin{tabular}{|c|l|c|}
        \hline    Step $\D$ & $98\%-$SD-Error$(\theta=0)$ & Order \\
      \hline              $2^{-5}$  & $0.034005669022654$ &  \\
            \cline{1-2}   $2^{-7}$  & $0.0344244107574687$ & $0.113$\\
              \cline{1-2} $2^{-9}$  & $0.0342415044537196$ & $(0.047)$\\
              \cline{1-2} $2^{-11}$ & $0.0331273071237492$ &\\
           \cline{1-2}   $2^{-13}$ & $0.0250195459354763$&\\
 \hline
  \end{tabular}
    \caption{\small The performance of fully explicit SD scheme (\ref{eq74}) applied to SDE (\ref{eq2}) with parameter set $(x_0, k_1, k_2, k_3, q, T)=(\frac{1}{16}, \frac{1}{16}, 1, 0.4, \frac{3}{4}, 1)$} \label{tab:exp_SD_0.75}
\end{table}


\item In practice, the computer time consumed to provide a desired level of accuracy, is of great importance. Especially, in financial applications, a scheme is considered better when except of its accuracy, it is implemented faster. As mentioned before, the SD method as well as the HAL method performs well in that aspect, compared to the implicit ALF method, which requires the estimation of a root of a nonlinear equation in each step and is therefore time consuming.  This is presented in  Table \ref{tab:order_time_SD_0.75} and Figure \ref{CEVModel_075wrtTIME}
which illustrates the advantage of the semi-discrete method SD, performing slightly better than HAL and BMM, better than BIM, and of course a lot better compared with ALF (over $1000$ times quicker to achieve an accuracy of almost $2$ decimal digits.)  Moreover, the explicit SD, performs slightly better in that aspect, as shown in Table \ref{tab:order_time_expSD_0.75}.
    \begin{table}[htbp]
    \centering
        \begin{tabular}{|c|c|c|}
        \hline          $q=0.75$ & Step $\D$ & Fully Explicit SD (Implicit)\\
    \hline         Time/Path(in sec)  & $2^{-5}$ &  $0.000013\quad  (0.000013)$ \\
      \hline  Time/Path(in sec)  & $2^{-7}$ &  $0.0000411\quad (0.0000422)$ \\
       \hline Time/Path(in sec)  & $2^{-9}$ &  $0.0001545\quad (0.0001586)$\\
  \hline      Time/Path(in sec)  & $2^{-11}$ & $0.0006048\quad (0.0006243)$\\
 \hline       Time/Path(in sec)  & $2^{-13}$ & $0.0024319\quad (0.0024975)$ \\
    \hline
  \end{tabular}
    \caption{\small Average computational time for a path (in seconds) for fully explicit SD method for $q=0.75.$}
    \label{tab:order_time_expSD_0.75}
\end{table}

\item A negative step of a numerical method appears when the computer-generated random variable exceeds a certain threshold, which tends to increase as the step size $\D$ decreases.
 Thus, the undesirable effect of negative values that are produced by some numerical schemes (such as the explicit Euler (EM) and standard Milsten (M) ), tends to disappear, since after a certain small step size, the threshold exceeds the maximum standard normal random number  attainable by the computer system.  
\end{itemize}

\section{Approximation of stochastic model (\ref{eq1}).}\label{sec:s7}
\setcounter{equation}{0}

So far we have focused on the process $(V_t),$ which is one part of the two-dimensional system
(\ref{eq1}). Nevertheless, it can be treated independently, since the only way that it interacts with the process $(S_t)$ is through the correlation $\rho$ of the Wiener processes.
 First we apply Ito's formula on  $\ln (S_t)$ to get,
\beqq\label{eq82}
\ln S_t=\ln S_0 + \int_{0}^{t}\mu du - \frac{1}{2}\int_{0}^{t}(V_u)^{2p}du + \int_{0}^{t}(V_u)^{p}dW_u , \quad t\in [0,T].
\eeqq
 
Then, we consider two different schemes for the integration of (\ref{eq82}).\footnote{The reason for not considering other schemes such as the two-dimensional Milstein is that they generally are time consuming, since they involve additional random number generation for the approximation of double Wiener integrals.} The first is the EM scheme which reads
\beqq\label{eq83}
\ln S_{t_{n+1}}^{EM}=\ln S_{t_n} + \mu \D - \frac{1}{2} (V_{t_n})^{2p} \D  +  (V_{t_n})^{p}\D W_n,
\eeqq
has strong convergence order $1/2$ and is easy to implement. The second scheme, which is based on an interpolation of the drift term and an interpolation of the diffusion term, considering decorrelation of the diffusion term, including a higher order Milstein term \cite[Sec.4.2]{kahl_jackel:2006}, is denoted IJK and is given by  \cite[Rel.(137)]{kahl_jackel:2006}

\beam\nonumber
\ln S_{t_{n+1}}^{IJK}&=&\ln S_{t_n} + \mu \D - \frac{1}{4} \left( (V_{t_n})^{2p} + (V_{t_{n+1}})^{2p} \right) \D  +  \rho(V_{t_n})^{p}\D \wt{W}_n\\
\label{eq84}&&  + 
\frac{1}{2} \left( (V_{t_n})^{p} + (V_{t_{n+1}})^{p} \right)(\D W_n - \rho\D \wt{W}_n)
 +\frac{1}{2} \rho pk_3(V_{t_n})^{q+p-1} \left( (\D \wt{W}_n)^2  - \D\right).
\eeam

\begin{figure}[ht]
  \caption{ \small Strong convergence error of the financial underlying process  $(S_t),$ as a function of CPU time (in sec) using log-Euler or IJK method with SD or BMM scheme for (\ref{eq1}) with $(x_0,S_0, \mu,k_1, k_2, k_3, q, T)=(\frac{1}{16}, 100, 0.05, \frac{1}{16}, 1, 0.4, \frac{3}{4}, 1),$ correlation $\rho=0$  and  $32$ digits of accuracy.}
  \centering
   \includegraphics[width=0.5\textwidth]{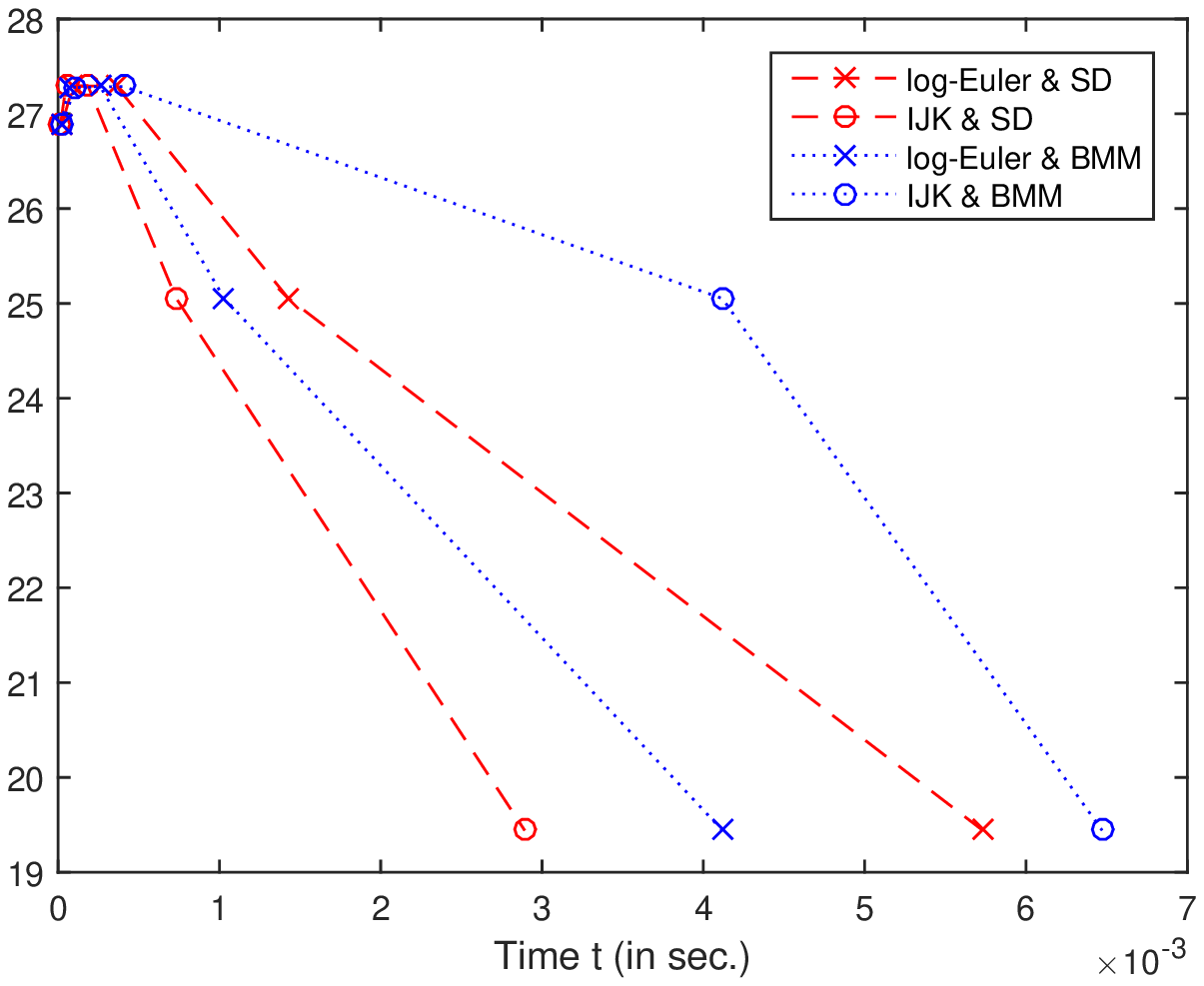}\label{StochVolModel_00wrtTIME}
  \end{figure}

\begin{figure}[ht]
  \caption{ \small Strong convergence error of the financial underlying process  $(S_t),$ as a function of CPU time (in sec) using log-Euler or IJK method with SD or BMM scheme for (\ref{eq1}) with $(x_0,S_0, \mu,k_1, k_2, k_3, q, T)=(\frac{1}{16}, 100, 0.05, \frac{1}{16}, 1, 0.4, \frac{3}{4}, 1),$ correlation $\rho=-0.4$  and  $32$ digits of accuracy.}
  \centering
   \includegraphics[width=0.5\textwidth]{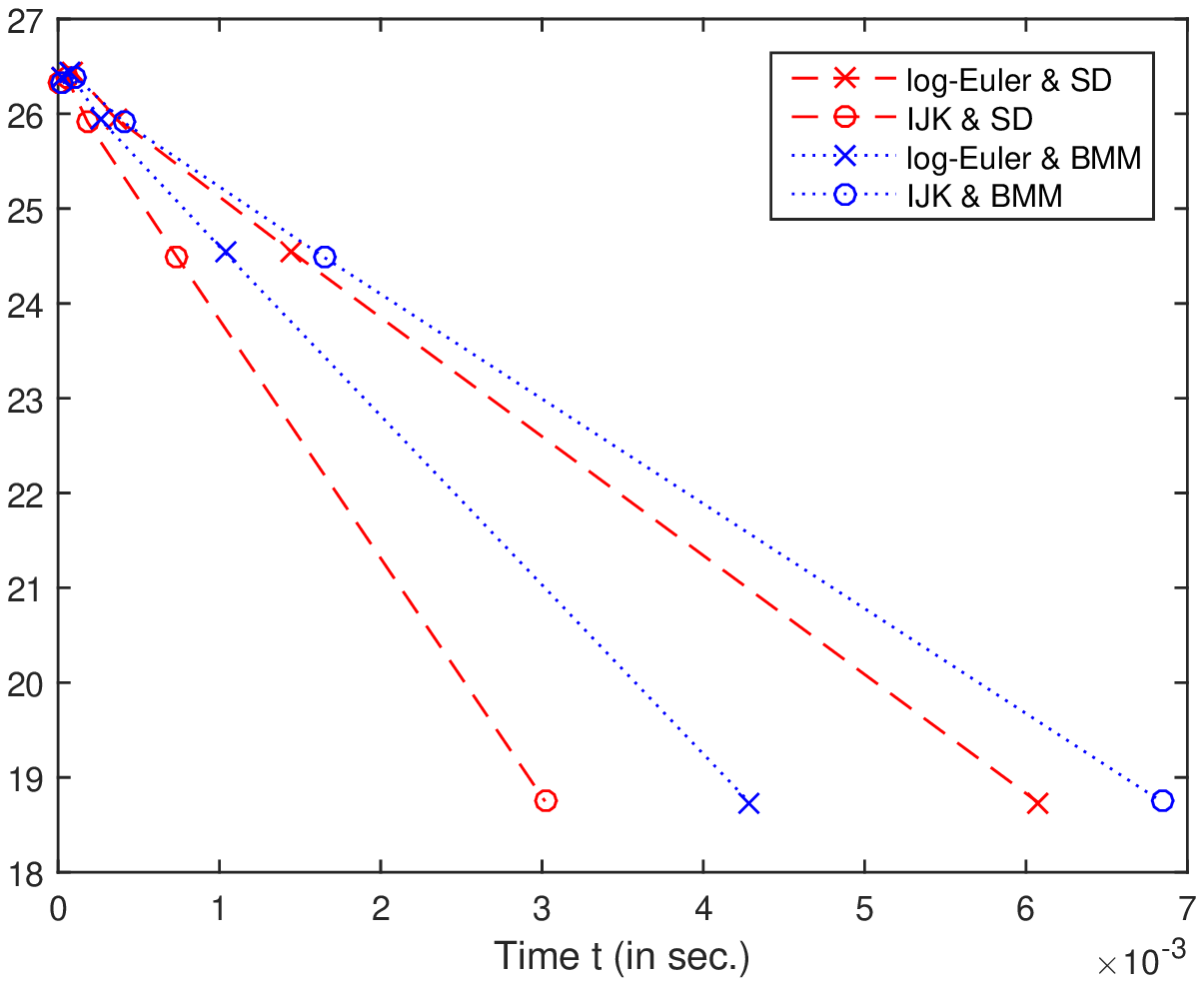}\label{StochVolModel_04wrtTIME}
  \end{figure}

\begin{figure}[ht]
  \caption{ \small Strong convergence error of the financial underlying process  $(S_t),$ as a function of CPU time (in sec) using log-Euler or IJK method with SD or BMM scheme for (\ref{eq1}) with $(x_0,S_0, \mu,k_1, k_2, k_3, q, T)=(\frac{1}{16}, 100, 0.05, \frac{1}{16}, 1, 0.4, \frac{3}{4}, 1),$ correlation $\rho=-0.8$  and  $32$ digits of accuracy.}
  \centering
   \includegraphics[width=0.5\textwidth]{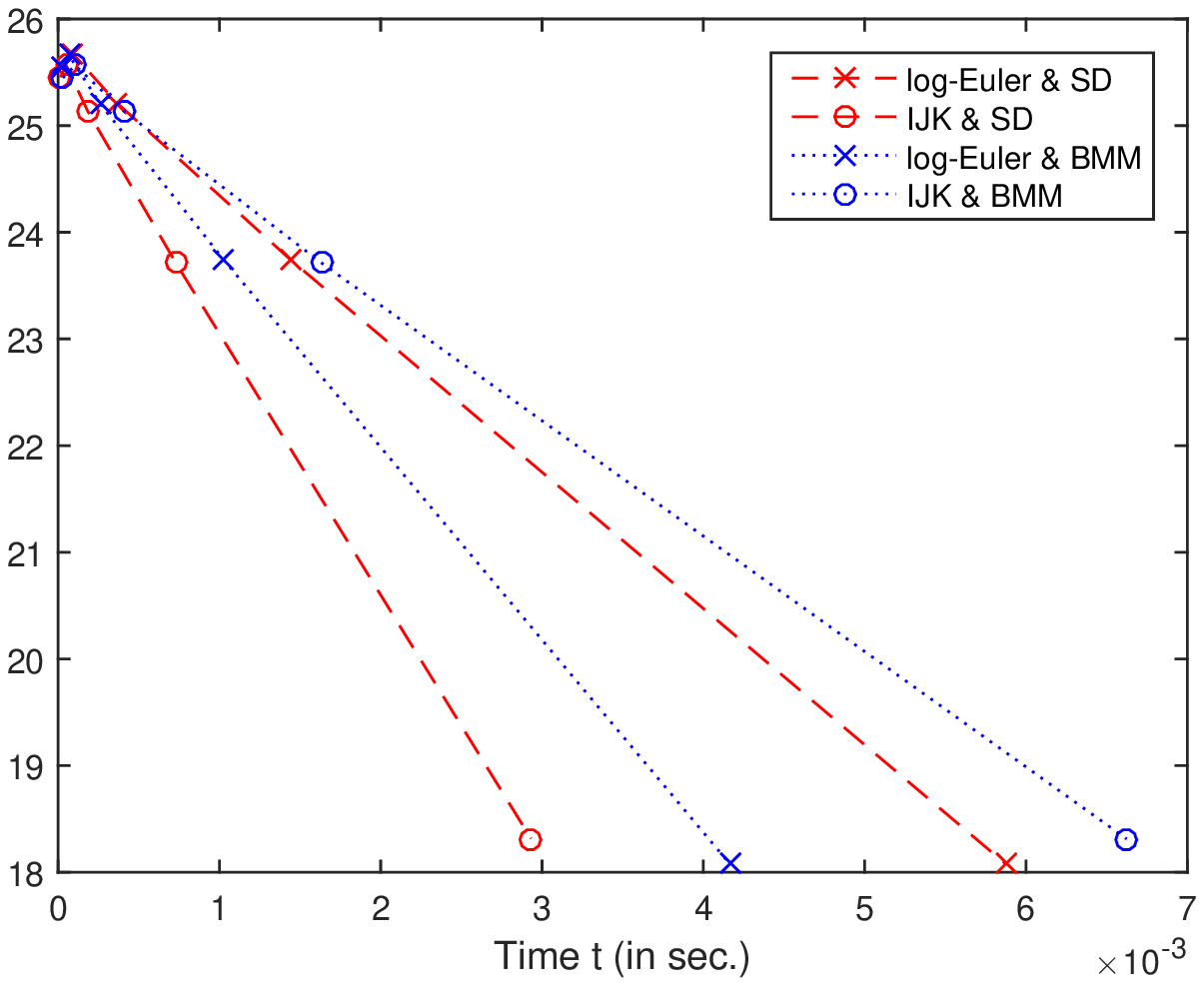}\label{StochVolModel_08wrtTIME}
  \end{figure}

We therefore consider the EM scheme (\ref{eq83}) combined with SD (\ref{eq74}), the IJK scheme (\ref{eq84}) combined with SD (\ref{eq74}) and compare with the case where the stochastic variance $(p=\frac{1}{2})$ is integrated with BMM scheme (\ref{eq81}), for three different correlation parameters, $\rho=0, \rho= -0.4$ and $\rho=-0.8$ with $S_0=100,\mu =0.05,$ as in \cite[Sec.5]{kahl_jackel:2006}. We present in Tables \ref{tab:2system0errors_0.75}, \ref{tab:2system0.4errors_0.75} and \ref{tab:2system0.8errors_0.75} and Figures \ref{StochVolModel_00wrtTIME}, \ref{StochVolModel_04wrtTIME} and \ref{StochVolModel_08wrtTIME}, the errors, in the sense of distance (\ref{l2dist}), for all the above considered ways of 
numerical integration of process $(S_t),$ for different step sizes, as well as the average computational time (in seconds) consumed for each discretization.

\begin{table}[htbp]  \scriptsize
\centering
        \begin{tabular}{|c|l|l|l|l|}
        \hline  Step $\D$ & EM$\&$SD-Error$(\theta=0.5)$ & IJK$\&$SD-Error$(\theta=0.5)$ & EM$\&$BMM-Error$(\Theta=0.5)$ & IJK$\&$ BMM-Error$(\Theta=0.5)$\\
 \hline $2^{-5}$  & $26.901\,(0.0000261)$& $26.901\,(0.0000159)$& $26.891\,(0.00002)$ &$26.890\,(0.0000294)$\\
\hline $2^{-7}$   & $27.288\,(0.0000919)$& $27.288\,(0.0000492)$& $27.277\,(0.0000676)$ &$27.277\,(0.0001043)$\\
\hline $2^{-9}$   & $27.298\,(0.0003595)$& $27.297\,(0.0001843)$& $27.289\,(0.0002610)$ &$27.288\,(0.0004081)$\\
\hline $2^{-11}$  & $25.057\,(0.0014255)$& $25.058\,(0.0007309)$& $25.051\,(0.0010309)$ &$25.051\,(0.0016191)$\\
\hline $2^{-13}$  & $19.441\,(0.0057322)$& $19.441\,(0.0028928)$& $19.442\,(0.0041177)$ &$19.442\,(0.0064721)$\\
 \hline  \end{tabular}
    \caption{\small $98\%-$Error, step size and average computational time of numerical integration of process $(S_t)$ using log-Euler or IJK method with SD or BMM scheme for (\ref{eq1}) with $(x_0,S_0, \mu,k_1, k_2, k_3, q, T)=(\frac{1}{16}, 100, 0.05, \frac{1}{16}, 1, 0.4, \frac{3}{4}, 1),$ correlation $\rho=0$  and  $32$ digits of accuracy.}
    \label{tab:2system0errors_0.75}
\end{table}

\begin{table}[htbp]  \scriptsize
\centering
        \begin{tabular}{|c|l|l|l|l|}
        \hline  Step $\D$ & EM$\&$SD-Error$(\theta=0.5)$ & IJK$\&$SD-Error$(\theta=0.5)$ & EM$\&$BMM-Error$(\Theta=0.5)$ & IJK$\&$ BMM-Error$(\Theta=0.5)$\\
 \hline $2^{-5}$  & $26.382\,(0.0000266)$ & $26.331\,(0.0000161)$ & $26.372\,(0.0000202)$ & $26.324\,(0.00003)$\\
\hline $2^{-7}$   & $26.448\,(0.0000951)$ &  $26.396\,(0.000005)$& $26.439\,(0.0000691)$ & $26.389\,(0.0001081)$\\
\hline $2^{-9}$   & $25.951\,(0.0003631)$ &  $25.909\,(0.000184)$& $25.944\,(0.0002606)$ & $25.904\,(0.0004131)$\\
\hline $2^{-11}$  & $24.540\,(0.0014506)$ & $24.494\,(0.0007355)$& $24.531\,(0.0010378)$  & $24.486\,(0.0016495)$\\
\hline $2^{-13}$  & $18.738\,(0.0060748)$ & $18.749\,(0.0030185)$ &$18.735\,(0.0042868)$  & $18.747\,(0.0068395)$\\
 \hline  \end{tabular}
    \caption{\small $98\%-$Error, step size and average computational time of numerical integration of process $(S_t)$ using log-Euler or IJK method with SD or BMM scheme for (\ref{eq1}) with $(x_0,S_0, \mu,k_1, k_2, k_3, q, T)=(\frac{1}{16}, 100, 0.05, \frac{1}{16}, 1, 0.4, \frac{3}{4}, 1),$ correlation $\rho=-0.4$  and  $32$ digits of accuracy.}
    \label{tab:2system0.4errors_0.75}
\end{table}

\begin{table}[htbp]  \scriptsize
\centering
        \begin{tabular}{|c|l|l|l|l|}
        \hline  Step $\D$ & EM$\&$SD-Error$(\theta=0.5)$ & IJK$\&$SD-Error$(\theta=0.5)$ & EM$\&$BMM-Error$(\Theta=0.5)$ & IJK$\&$ BMM-Error$(\Theta=0.5)$\\
 \hline $2^{-5}$  & $25.552\,(0.0000263)$ &$25.455\,(0.0000159)$ &$25.541\,(0.0000199)$ &$25.449\,(0.0000296)$\\
\hline $2^{-7}$   & $25.670\,(0.0000932)$ &$25.569\,(0.0000494)$ &$25.659\,(0.0000678)$ &$25.564\,(0.0001059)$\\
\hline $2^{-9}$   & $25.217\,(0.0003622)$ &$25.137\,(0.0001835)$ &$25.208\,(0.0002595)$ &$25.132\,(0.0004111)$\\
\hline $2^{-11}$  & $23.743\,(0.0014407)$ &$23.711\,(0.0007306)$ &$23.734\,(0.0010307)$ &$23.707\,(0.0016376)$\\
\hline $2^{-13}$  & $18.082\,(0.005871)$  &$18.316\,(0.0029312)$ &$18.078\,(0.0041637)$ &$18.312\,(0.0066239)$\\
 \hline  \end{tabular}
    \caption{\small $98\%-$Error, step size and average computational time of numerical integration of process $(S_t)$ using log-Euler or IJK method with SD or BMM scheme for (\ref{eq1}) with $(x_0,S_0, \mu,k_1, k_2, k_3, q, T)=(\frac{1}{16}, 100, 0.05, \frac{1}{16}, 1, 0.4, \frac{3}{4}, 1),$ correlation $\rho=-0.8$  and  $32$ digits of accuracy.}
    \label{tab:2system0.8errors_0.75}
\end{table}
Figures 	\ref{StochVolModel_00wrtTIME}, \ref{StochVolModel_04wrtTIME} and \ref{StochVolModel_08wrtTIME} indicate that in all cases the favorable choice is to integrate $(S_t)$ using IJK method combined with the SD scheme for $(V_t)$ in model (\ref{eq1}).
 The above combination seems to be the better one, w.r.t. CPU time, for every correlation coefficient considered. 	

\section{Conclusion.}\label{sec:s8}
\setcounter{equation}{0}

In this paper, we exploit further the semi-discrete method (SD), originally appeared in \cite{halidias:2012}, to numerically approximate stochastic processes that appear in financial mathematics and are meant to be non negative. In \cite{halidias_stamatiou:2015} we examined the Heston $3/2-$model, that is a mean reverting process with super-linear diffusion, described by a SDE of the form (\ref{eq2}) with $q=3/2.$ Now, we deal with SDEs with sub-linear diffusion coefficients of the type $(x_t)^q$ with $1/2<q<1.$ These kind of SDEs, called mean reverting CEV processes, appear in stochastic models, where they represent the instantaneous volatility-variance of an underlying financially observable variable. We prove theoretically the strong convergence of our proposed SD scheme, revealing the order of convergence. The resulting polynomial rate in Theorem \ref{t3}, may not appear appealing at first sight, because of its low magnitude. Nevertheless, as it is shown in the numerical experiment section, where a comparative study is presented between various positivity preserving schemes, the SD method seems to be the best w.r.t. CPU time consumption. The advantage of the SD method here is that although implicit, has an explicit formula and thus requires fewer arithmetic operations and consequently less computational time. Moreover, our method can cover cases where (\ref{eq2}) has time varying coefficients, i.e. $k_1(t), k_2(t), k_3(t).$

We also treat the whole two-dimensional stochastic volatility model (\ref{eq1}). In order to do that, we actually integrate 
the process $\ln (S_t)$ which satisfies a SDE of the form (\ref{eq82}) and in the end transform back for $(S_t).$
 We only consider two different schemes for the integration of $\ln (S_t),$ namely the Euler Maruyama (EM) scheme, which is easy to implement and the IJK scheme \cite[Rel.(137)]{kahl_jackel:2006} which is shown to be the most efficient method, robust and simple as EM \cite{kahl_jackel:2006}. We do not apply other two-dimensional schemes, such as for example the  Milstein scheme, since they are in general time consuming, as they involve approximations of double Wiener integrals which require
additional random number generation.  We therefore combine the EM scheme with SD ((\ref{eq83}) $\&$ (\ref{eq74})), the IJK scheme with SD ((\ref{eq84}) $\&$ (\ref{eq74})) and compare with the case where the stochastic variance $(p=\frac{1}{2})$ is integrated with BMM scheme (\ref{eq81}), for three different correlation parameters, $\rho=0, \rho= -0.4$ and $\rho=-0.8$ with $S_0=100,\mu =0.05,$ as in \cite[Sec.5]{kahl_jackel:2006}.  The combination IJK with SD seems to be the most favorable w.r.t. CPU time, for all the cases. 	


\appendix
\section{Some numerical schemes for the integration of the variance-volatility process $(V_t).$}\label{A1}
\setcounter{equation}{0}

We consider a partition of the time interval $[0,T]$ with $0=t_0<t_1<...<t_N=T$ and discretization steps $\D_n:=t_{n+1}-t_n$ for $n=0,\ldots,N-1.$ Moreover, we denote by $\D W_n:=W_{t_{n+1}}- W_{t_{n}}$ the increments of the Brownian motion. We show in the following subsections some numerical schemes for the approximation of
\beqq\label{vol}
V_t = V_0 + \int_0^t (k_1-k_2 V_s)ds + \int_{0}^{t}k_3(V_s)^q dW_s, \quad t\in [0,T]
\eeqq
and make some brief comments on them. We also denote $V_n:=V_{t_n}.$

\subsection*{Standard Euler-Maruyama scheme}\label{A1.1}

The Euler method, applied to the SDE setting, already appeared in the $50's$ through Maruyama \cite{maruyama:1955} and thereafter there has been an extensive study on numerical approximations of solutions of SDEs (we just mention \cite{kloeden_neuenkirch:2013} for a recent review on numerical methods for SDEs with applications in finance and references therein).

The explicit Euler-Maruyama (EM) scheme for the process $(V_t)$ is given by 
\beqq\label{EM}
V_{n+1}^{EM}=V_n +(k_1-k_2 V_n)\D_n + k_3(V_n)^q \D W_n,
\eeqq
for $n=0,\ldots,N-1.$ Clearly $\bfP(V_{n+1}<0 | V_n>0)>0,$ thus the EM scheme can produce negative values with positive probability, or in the notion of 
 \cite{schurz:1996} we say that (\ref{EM}) has a finite life time.

\subsection*{Standard Milstein scheme}\label{A1.2}

The standard one dimensional Milstein (M) scheme 
contains some extra terms derived by Ito-Taylor expansion \cite[Sec.5]{kloeden_platen:1995}, and applied to $(V_t)$ reads  
\beqq\label{M}
V_{n+1}^{M}=V_n +(k_1-k_2 V_n)\D_n + k_3(V_n)^q \D W_n + \frac{1}{2}(k_3)^2q(V_n)^{2q-1}\left((\D W_n)^2-\D_n\right),
\eeqq
for $n=0,\ldots,N-1$ where we have retained terms of order $(\D_n).$ Again (M) scheme has a finite life time.

\subsection*{Balanced Implicit Method}\label{A1.3}

The balanced implicit method (BIM) \cite[Rel (3.2)]{milstein_et_al:1998} was the first attempt to treat the problem of invariance-preserving of specific domains of the underlying process and reads
\beqq\label{BIM}
V_{n+1}^{BIM}=V_n +(k_1-k_2 V_n)\D_n + k_3(V_n)^q \D W_n + \left(c^0(V_n)\D_n + c^1(V_n)|\D W_n|\right)(V_n - V_{n+1}),
\eeqq
for $n=0,\ldots,N-1$ where $c^0$ and $c^1$ are appropriate weight functions. The choice $c^0(x)=k_2$ and $c^1(x)=k_3x^{q-1}$ preserves positivity \cite[Sec. 5]{kahl_schurz:2006}. Rearranging the above equation, we get the expression
\beqq\label{BIM2}
V_{n+1}^{BIM}=\frac{V_n + k_1\D_n + k_3(V_n)^q (\D W_n + |\D W_n|) }{1 +k_2\D_n + k_3(V_n)^{q-1} |\D W_n|}.
\eeqq

\subsection*{Balanced Milstein Method}\label{A1.4}

The balanced Milstein method (BMM), was proposed in \cite{kahl_schurz:2006}, for an improvement of the BIM  in the stability behavior but also both in the rate of convergence.
It is given by the following linear implicit relation
\beao
V_{n+1}^{BMM}&=&V_n +(k_1-k_2 V_n)\D_n + k_3(V_n)^q \D W_n + \frac{1}{2}(k_3)^2q(V_n)^{2q-1}\left( (\D W_n)^2 - \D_n\right)\\
&& +  \left(d^0(V_n)\D_n + d^1(V_n)( (\D W_n)^2 - \D_n)\right)(V_n - V_{n+1}),
\eeao
for $n=0,\ldots,N-1$ where $d^0$ and $d^1$ are appropriate weight functions. The choice $d^0(x)=\Theta k_2 + \frac{1}{2}(k_3)^2q|x|^{2q-2},$ where $\Theta\in[0,1]$ and $d^1(x)=0$ implies an eternal life time for the scheme \cite[Th. 5.9]{kahl_schurz:2006}, in the sense that $\bfP(V_{n+1}>0 | V_n>0)=1.$ The step sizes $\D_n$ have to be such that $\D_n < \frac{2q-1}{2qk_2(1-\Theta)}.$ The relaxation parameter resembles to the implicitness parameter ($ \theta$ in our notation). For $\Theta=1$ there is no restriction in the step size, but it is recommended when possible \cite[Rem. 5.10]{kahl_schurz:2006} to take $\Theta=1/2.$ 
 Rearranging  with the above specifications leads to
\beqq\label{BMM}
V_{n+1}^{BMM}=\frac{V_n +(k_1-(1-\Theta)k_2 V_n)\D_n + k_3(V_n)^q \D W_n + \frac{1}{2}(k_3)^2q(V_n)^{2q-1}  (\D W_n)^2}{1 + \Theta k_2\D_n + \frac{1}{2}(k_3)^2q|V_n|^{2q-2}\D_n}.
\eeqq

Finally, the proposed semi-discrete (SD) scheme reads 
\beqq\label{SD}
V_{n+1}^{SD}=\left( \sqrt{V_{n}\left(1- \frac{k_2\D}{1+k_2\theta\D}\right) + \frac{k_1\D}{1+k_2\theta\D} - \frac{(k_3)^2\D}{4(1+k_2\theta\D)^2}(V_{n})^{2q-1}} +\frac{k_3}{2(1+k_2\theta\D)}(V_{n})^{q -\frac{1}{2}}\D W_n\right)^2.
\eeqq

Increasing the time horizon $T$ results in an increase of the percentage of negative paths of EM and M. On the other hand BIM, BMM and of course SD are not affected by that, since they preserve their positivity on any interval $[0,T].$


\begin{thebibliography}{99}
\baselineskip12pt


\bibitem{alfonsi:2013}
{\sc A. Alfonsi,}
{\em Strong order one convergence of a drift implicit Euler scheme: Application to  the CIR process,}
Stat. Prob. Let., 83\ (2013),
pp. 602-607.

\bibitem{andersen_piterbarg:2007}
{\sc L.B.G. Andersen, V.V. Piterbarg,}
{\em Moment explosions in stochastic volatility models,}
 Finance Stoch., 11\ (2007),
pp. 29-50.


\bibitem{berkaoui:2004}
{\sc A. Berkaoui,}
{\em Euler scheme for solutions of stochastic differential equations,}
Portugalia Mathematica Journal, 61,\ (2004),
pp. 461-478.

\bibitem{brennan_schwartz:1980}
{\sc M.J. Brennan, E.S. Schwartz,}
{\em Analyzing convertible bonds,}
Journal of Financial and Quantitative Finance, 4,\ (1980),
pp. 907-929.


\bibitem{chan_et_al:1992}
{\sc K.C. Chan, G.A Karolyi, F.A. Longstaff, A.B. Sanders,}
{\em An empirical comparison of short-term interest rate,}
 Journal of Finance, 47(3), \ (1992),
pp. 1209-1227.

\bibitem{cox_et_al:1985}
{\sc J.C. Cox, J.E. Ingersoll, S.A. Ross,}
{\em A theory of the term structure of interest rates,}
 Econometrica, 53, \ (1985),
pp. 385-407.
  
\bibitem{cox_et_al:2013}
{\sc S.G. Cox, M. Hutzenthaler, A. Jentzen,}
{\em Local Lipschitz continuity in the initial value and strong completeness for nonlinear stochastic differential equations,}
 arXiv:1309.5595v1, \ (2013).


\bibitem{dereich_et_al:2011}
{\sc S. Dereich, A. Neunkirch, L. Szpruch,}
{\em An Euler-type method for the strong approximation of the Cox-Ingersoll-Ross process,}
Proceedings of The Royal Society, \ (2011).
%

\bibitem{gronwall:1919}
{\sc T.H. Gronwall,}
{\em Note on the derivatives with respect to a parameter of the solutions of a system of differential equations,}
 Annals of Mathematics, 20, \ (1919),
pp. 292-296.

\bibitem{halidias:2012}
{\sc N. Halidias,}
{\em Semi-discrete approximations for stochastic differential equations and applications,}
 International Journal of Computer Mathematics, 89(6), \ (2012),
pp. 1-15.

\bibitem{halidias:2015a}
{\sc N. Halidias,}
{\em A new numerical scheme for the CIR process,}
 Monte Carlo Methods Appl., \ (2015a).

\bibitem{halidias:2015b}
{\sc N. Halidias,}
{\em An explicit and positivity preserving numerical scheme for the mean reverting CEV model,}
 http://arxiv.org/pdf/1501.03434, \ (2015b).

\bibitem{halidias_stamatiou:2015}
{\sc N. Halidias, I.S. Stamatiou,}
{\em On the numerical solution of some nonlinear stochastic differential equations using the semi-discrete method,}
to appear in Computational Methods in Applied Mathematics, Special Issue, \ (2015).
 
\bibitem{hurd_kuznetsov:2008}
{\sc T.R. Hurd, A. Kuznetsov,}
{\em Explicit formulas for Laplace transforms of stochastic integrals,}
 Markov Process. Relat. Fields, 14, \ (2008),
pp. 277-290.


\bibitem{kahl_jackel:2006}
{\sc C. Kahl, P. Jackel,}
{\em Fast strong approximation Monte Carlo schemes for stochastic volatility models,}
 Quantitative Finance, 6, \ (2006),
pp. 513-536.

\bibitem{kahl_schurz:2006}
{\sc C. Kahl, H. Schurz,}
{\em Balanced Milstein Methods for ordinary SDEs,}
 Monte Carlo Methods and Appl, 12(6), \ (2006),
pp. 143-170.
 
\bibitem{karatzas_shreve:1988}
{\sc I. Karatzas, S.E. Shreve,}
{\em Brownian motion and stochastic calculus,}
Springer-Verlag New York, \ (1988).
%
\bibitem{kloeden_neuenkirch:2013}
{\sc P. Kloeden, A. Neuenkirch,}
{\em Convergence of numerical methods for stochastic differential equations in mathematical finance,}
 Recent developments in Computational Finance, (T. Gerstner and P. Kloeden, eds), \ (2013),
pp. 49-80.

\bibitem{kloeden_platen:1995}
{\sc P. Kloeden, E. Platen,}
{\em Numerical solution of stochastic differential equations,}
Vol 23, Stochastic Modeling and Applied Probability, Springer-Verlag Berlin, corrected 2nd printing, \ (1995).

\bibitem{kloeden_platen_schurz:2003}
{\sc P. Kloeden, E. Platen, H. Schurz,}
{\em Numerical solution of stochastic differential equations through computer experiments,}
Springer-Verlag Berlin, corrected 3rd printing, \ (2003).

\bibitem{mao:1997}
{\sc X. Mao,}
{\em Stochastic Differential Equations and Applications,}
Horwood Publishing, \ (1997).

\bibitem{maruyama:1955}
{\sc  G. Maruyama,}
{\em Continuous Markov processes and stochastic equations,}
 Rend. Circ. Mat. Palermo, 4(1),\ (1955),
pp. 48-90.

\bibitem{milstein_et_al:1998}
{\sc G.N. Milstein, E. Platen, H. Schurz,}
{\em Balanced implicit methods for stiff stochastic systems,}
SIAM J. Numer. Anal., 35(3),\ (1998),
pp. 1010-1019.

\bibitem{olver:1997}
{\sc F.W.J. Olver,}
{\em Asymptotics and special functions,}
AKP classics, Wellesley, Mass,\ (1997).



\bibitem{schurz:1996}
{\sc H. Schurz,}
{\em Numerical regularization for SDEs: construction of nonnegative solutions,}
 Dyn. Systems Appl., 5, \ (1996), 
pp. 323-352.



\bibitem{tretyakov_zhang:2013}
{\sc M.V. Tretyakov, Z. Zhang,}
{\em A fundamental mean-square convergence theorem for SDEs with locally Lipshcitz coefficients and its applications,}
SIAM J. Numer. Anal., \ (2013),
pp. 3135-3162.

\bibitem{yamada_watanabe:1971}
{\sc T. Yamada, S. Watanabe,}
{\em On the uniqueness of solutions of stochastic differential equations,}
 J. Math. Kyoto Univ., 11, \ (1971), pp. 155-167.
\end{thebibliography}
\end{document}